\documentclass[absolute]{mymathart}
\usepackage[theorems]{mymathmacros}
\usepackage[usenames,dvipsnames]{color}
\usepackage{subfig}
\usepackage[shortlabels]{enumitem}
\usepackage{hyperref}
\usepackage{amsmath}
\usepackage{amssymb}
\usepackage{graphicx}
\usepackage{epigraph}
\usepackage{amsfonts,amssymb,color}
\usepackage[mathscr]{eucal}
\usepackage{amsmath, amsthm}
\usepackage{mathrsfs}
\usepackage{amsbsy}
\usepackage{float}
\usepackage{verbatim}
\usepackage{amsfonts} 
\usepackage{pifont}

\title[Rigidity of pressures of H\"older potentials]{Rigidity of pressures of H\"older potentials and the fitting of analytic functions via them}

\author{Liangang Ma \lowercase{and} Mark Pollicott}

\address{Liangang Ma, School of Mathematics and Statistics, Ludong University, Yantai 264025, Shandong, P. R. China.} 
\email{maliangang000@163.com}

\address{Mark Pollicott, Department of Mathematics, Warwick University, Coventry, CV4 7AL, UK.} 
\email{masdbl@warwick.ac.uk}

\thanks{The first author is supported by 12001056 from NSFC and ZR2019QA003 from SPNSF. The second author is partly supported by ERC-Advanced Grant 833802-Resonances and EPSRC Grant EP/T001674/1.}  

\newtheorem{theorem}[subsection]{Theorem}
\newtheorem{lemma}[subsection]{Lemma}
\newtheorem{pro}[subsection]{Proposition}
\newtheorem{prob}[subsection]{Problem}
\newtheorem*{CLT}{Central Limit Theorem}
\newtheorem{coro}[subsection]{Corollary}

\newtheorem{rem}[subsection]{Remark}

\numberwithin{equation}{section}

\begin{document} 

\begin{abstract}
The first part of this work is devoted to the study of higher differentials of pressure functions of H\"older potentials on shift spaces of finite type. By describing the differentials of pressure functions via the Central Limit Theorem for the associated random processes, we discover some rigid relationships between differentials of various orders. The rigidity imposes obstructions on fitting candidate convex analytic functions by pressure functions of H\"older potentials globally, which answers a question of Kucherenko-Quas. In the second part of the work we consider fitting candidate analytic germs by pressure functions of locally constant potentials. We prove that all 1-level candidate germs can be realised by pressures of some locally constant potentials, as long as number of the symbolic set is large enough. There are also some results on fitting 2-level germs by pressures of locally constant potentials obtained in the work.         
\end{abstract}
 
 \maketitle

\section{Introduction}\label{sec1}

This work deals with traditional topics in thermodynamic formalism \cite{Bow, Rue1}, which originates from theoretical physics. We focus on shift spaces of finite type here, which model dynamics of  some smooth systems such as Axiom-A Diffeomorphisms   through Markov partitions. Given a symbolic set $\Lambda$ of finite symbols and a continuous potential (observable) $\phi$ on the shift space $\Lambda^\mathbb{N}$, a core concept in thermodynamic formalism is the pressure $P(\phi)$. People are particularly interested in the pressure function $P(t\phi)$ with the variable $t>0$ representing the inverse temperature.  A sharp change in the pressure function (or other terms) is usually called a phase transition as $t$ varies, see for example \cite{IRV, IT1, IT2, KQW, Lop1, Lop2, Lop3, Sar}.  

For H\"older continuous potentials, Ruelle \cite{Rue2} proved that the pressure function $P(t\phi)$ is analytic for $t\in(0,\infty)$ (in fact he proved that $P(\psi)$ depends analytically on $\psi$ for $\psi$ in the H\"older space $C^{0,h}(X)$ with $X$ being a transitive subshift space of finite type and $0<h\leq 1$ being the exponent \cite{GT}). A key ingredient in his proof is the use of Ruelle (transfer) operator \cite{BDL, GLP} acting on functions in the H\"older space. Moreover, the equilibrium measure of $t\phi$ for any $t>0$ and H\"older potential $\phi$ is always unique, so there are in fact no phase transitions in this case. Let 
\begin{center}
$P^{(n)}(t)=P^{(n)}(t\phi)=\cfrac{d^nP(t\phi)}{dt^n}$
\end{center}
be the $n$-th differential of the pressure function $P(t\phi)$ with respect to $t\in(0,\infty)$ for some fixed H\"older potential $\phi$. We also write 
\begin{center}
$P^{(1)}(t)=P'(t), P^{(2)}(t)=P''(t), P^{(3)}(t)=P'''(t),\cdots$ 
\end{center}
intermittently in the following. We discover that there is some rigid relationship between the differentials of the pressure function.

\begin{theorem}\label{thm8}
For a H\"older potential $\phi$ on a full shift space of finite type, let $P(t)=P(t\phi)$ be its pressure. Then there exists some positive number $M_\phi$ depending on $\phi$, such that 
\begin{equation}\label{eq50}
\sqrt{2\pi^3}(P^{(2)}(t))^{3/2}|P^{(3)}(t)|\leq 9|P^{(3)}(t)|+2|P^{(4)}(t)| + 3\sqrt{2\pi^3} M_\phi (P^{(2)}(t))^{5/2}
\end{equation}
for any $t>0$.

\end{theorem}

 A potential $\phi$ is said to be \emph{generic} (or we say it defines a non-lattice distribution, cf. \cite{CP, Fel, PP}), if for any normalised potential $\psi$, the spectral radius of the complex Ruelle operator $\mathcal{L}_{\psi+it\phi}$ is less than $1$ for any $t\neq 0$. For pressure functions of generic potentials, Theorem \ref{thm8} can be strengthened to the following result.    

\begin{theorem}\label{thm5}
For a generic H\"older potential $\phi$ on a full shift space of finite type, let $P(t)=P(t\phi)$ be its pressure. Then there exists some positive number $M_\phi$ depending on $\phi$, such that 
\begin{equation}\label{eq19}
|P^{(3)}(t)\big(1-\sqrt{2\pi}(P^{(2)}(t))^{3/2}\big)|\leq 3 M_\phi P^{(2)}(t)
\end{equation}
for any $t>0$.

\end{theorem}

This means the second differential of the pressure function of a generic H\"older potential imposes some global subtle restriction on its third differential. It would be interesting to try to interpret the meaning of $P''(t)=\frac{1}{\sqrt[3]{2\pi}}$ for the pressure function at individual parameters. Let $\sigma: \Lambda^\mathbb{N}\rightarrow\Lambda^\mathbb{N}$  denote the shift map. Both the proofs of Theorem \ref{thm8} and \ref{thm5} require use of the Ruelle operator and the Central Limit Theorem (CLT) for the process $\{f\circ\sigma^n\}_{n\in\mathbb{N}}$, with the latter one depending on a finer CLT in the generic case. Recall that there are some expressions on the higher differentials of the pressure function by Kotani and Sunada in \cite{KS1} for smooth systems, and we refer the readers to \cite{KS2} for a CLT for random walks on crystal lattices.

It is well-known that $P(t\phi)$ is convex and Lipschitz for continuous $\phi$, moreover, the supporting lines of its graph must intersect the vertical axis in a closed bounded interval in $[0,\infty)$. Kucherenko and Quas have shown that any such function can be realised by the pressure function of some continuous potential on some shift space of finite type \cite[Theorem 1]{KQ}, whose result fits into Katok's flexibility programme \cite{BKR}. However, the continuous potentials constructed in their work are not H\"older, so they ask the following question (their original problem is set in the multidimensional case). 

\begin{prob}[Kucherenko-Quas]\label{prob1}
Can a convex, Lipschitz analytic function with its supporting lines intersecting the vertical axis in a closed bounded interval in $[0,\infty)$ be realised by the pressure function of some  H\"older potential on some shift space of finite type? 
\end{prob}

Our following results are dedicated to an answer to their problem. We first point out that any convex, Lipschitz analytic function with its supporting lines intersecting the vertical axis in a closed bounded interval in $[0,\infty)$ can be approximated by sequences of pressure functions of locally constant potentials  on some shift space of finite type.

\begin{coro}\label{thm1}
Let $F(t)$ be a convex Lipschitz function on $(\alpha,\infty)$ for some $\alpha>0$ with Lipschitz constant $L>0$, such that its supporting lines intersect the vertical axis in $[\underline{\gamma},\overline{\gamma}]$ with $0\leq\underline{\gamma}\leq\overline{\gamma}<\infty$. Then there exists a sequence of locally constant potentials $\{\phi_n\}_{n=1}^\infty$ on some shift space of finite type, such that
\begin{equation}\label{eq1}
\lim_{n\rightarrow\infty} P(t\phi_n)=F(t)
\end{equation}
for any $t\in(\alpha,\infty)$.

\end{coro}

\begin{proof}
This is an instant corollary of Kucherenko-Quas' result.  Let 
\begin{center}
$\Lambda=\{0,1,\cdots,\lfloor e^{\overline{\gamma}}\rfloor\}\times\{\lfloor \underline{\gamma}\rfloor,\cdots,\lceil \overline{\gamma}\rceil\}\}\times\{\lfloor -L\rfloor,\cdots,\lceil L\rceil\}$,
\end{center}
in which $\lfloor\ \rfloor$ and $\lceil\ \rceil$ represent the floor and ceiling function respectively.
According to \cite[Theorem 1]{KQ}, there exists a continuous potential $\phi_F: \Lambda^\mathbb{Z}\rightarrow \mathbb{R}$, such that 
\begin{center}
$P(t\phi_F)=F(t)$
\end{center}
on $(\alpha,\infty)$. Now let
\begin{center}
$\phi_n(x)\doteq\phi_{n,-}(x)=\inf\{\phi_F(x): x\in[x_{-n}x_{-n+1}\cdots x_n]\}$
\end{center}
for any $x=\cdots x_{-(n+1)} x_{-n}\cdots x_n x_{n+1}\cdots\in \Lambda^\mathbb{Z}$ and $n\in\mathbb{N}$, in which $[x_{-n}x_{-n+1}\cdots x_n]$ means the corresponding cylinder set.  $\phi_n$ is a locally constant potential for any fixed $n$. Now fix $t\in (\alpha,\infty)$, by properties of the pressure function (see for example \cite[6.8]{Rue1}),
\begin{equation}\label{eq53}
|P(t\phi_n)-P(t\phi_F)|\leq |t|\parallel\phi_n-\phi_F\parallel_{\infty}.
\end{equation}
Since $\phi_F$ is continuous, this implies (\ref{eq1}). 

\end{proof}

One can see that in the above proof the increasing sequence of pressures $\{P(t\phi_{n,-})\}_{n\in\mathbb{N}}$ satisfies
\begin{center}
$P(t\phi_{n,-})\nearrow F(t)$
\end{center}
as $n\rightarrow\infty$  since $\{\phi_{n,-}\}_{n\in\mathbb{N}}$ is an increasing sequence tending to $\phi_F$ (see \cite[Theorem 9.7(ii)]{Wal1}). Alternatively, one can take 
\begin{center}
$\phi_{n,+}(x)=\sup\{\phi_F(x): x\in[x_{-n}x_{-n+1}\cdots x_n]\}$,
\end{center}
which results in a decreasing sequence of locally constant potentials approximating $\phi_F(x)$, or  
\begin{center}
$\phi_{n,\pm}(x)=\cfrac{\phi_{n,-}(x)+\phi_{n,+}(x)}{2}$,
\end{center}
which also results in a sequence of locally constant potentials approximating $\phi_F(x)$, while their pressure functions both approximate $F(t)$. See Corollary \ref{cor3} for an interpretation
of the result from another point of view.

\begin{rem}
The convergence in Corollary \ref{thm1} is uniform for $t$ in a bounded domain since $\Lambda^\mathbb{Z}$ is a compact metric space by considering (\ref{eq53}). 
\end{rem}

\begin{rem}
A locally constant potential is of course H\"older, so according to  Ruelle's result, the pressure functions $\{P(t\phi_{n,-})\}_{n\in\mathbb{N}}$ are all analytic.

\end{rem}

The following result confirms that some convex analytic functions cannot be fitted by the pressure of any H\"older potential on any shift space of finite type, which gives a negative answer to Problem \ref{prob1}.

\begin{theorem}\label{thm7}
For any $\alpha>0$, there exists a strictly convex analytic function $F(t)$ on $(\alpha,\infty)$, with its supporting lines intersecting the vertical axis in $[\underline{\gamma},\overline{\gamma}]\subset [0,\infty)$, such that there does not exist any H\"older potential $\phi$ on any shift space of finite type satisfying
\begin{center}
$P(t\phi)= F(t)$
\end{center}
on $(\alpha,\infty)$.
\end{theorem}

For an explicit example of convex analytic functions in Theorem \ref{thm7}, one can simply take 
\begin{center}
$F_{2,3,1}(t)=\cfrac{2t^2+3t+te^{-t^2}+e^{-t^2}}{t}$
\end{center} 
on $(\alpha,\infty)$ for any $\alpha>0$. See Proposition \ref{pro3} for a family of such examples. Thus one can see that there are in fact elementary functions which cannot be fitted by pressures of H\"older potentials on shift spaces of finite type. 

In the following we consider fitting convex analytic functions locally instead of globally, only by pressures of locally constant potentials on shift spaces of finite type.  Let 
\begin{center}
$\Lambda_n=\{1,2,\cdots,n\}$
\end{center}
be the symbolic set of $n$ symbols.
\begin{theorem}\label{thm3}
Let $t_*>0$ and $(a_0,a_1)\in\mathbb{R}^2$ satisfying
\begin{equation}\label{eq32}
\frac{a_0}{t_*}>a_1.
\end{equation}
Then for any $n\in\mathbb{N}$ large enough, there exist some $0\leq m_{t_*,a_0,a_1,n}< M_{t_*,a_0,a_1,n}<\infty$ depending on $t_*, a_0, a_1, n$, such that for any $a_2\in[m_{t_*,a_0,a_1,n}, M_{t_*,a_0,a_1,n}]$, there exists some sequence of reals $\{c_{i,n}\}_{i=1}^n$, such that the locally constant potential
\begin{center}
$\phi(x)=c_{x_0,n}$
\end{center}  
for $x=\cdots x_{-1}x_{0}x_{1}\cdots\in[x_0]$ on the full shift space $\Lambda_n^{\mathbb{Z}}$ satisfies
\begin{equation}\label{eq41}
P(t\phi)=a_0+a_1(t-t_*)+\cfrac{a_2}{2!}(t-t_*)^2+O((t-t_*)^3)
\end{equation}
on $[t_*-\delta_n,t_*+\delta_n]$ for some $\delta_n>0$.
\end{theorem}

This means we can fit some germs at $t_*$ up to level $2$ by pressures of some locally constant potentials when the number of symbols of the shift space is large enough. The values $\delta_n, \{c_{i,n}\}_{i=1}^n$ all depend on $t_*, a_0, a_1, n$ and $a_2$ in fact, while we only indicate the dependence of $m_{t_*,a_0,a_1,n}$ and $M_{t_*,a_0,a_1,n}$ as we are particularly interested in their values in the context of Theorem \ref{thm3}. There are some results on the values of 
\begin{center}
$\{m_{t_*,a_0,a_1,n}, M_{t_*,a_0,a_1,n}\}_{n\in\mathbb{N}}$ 
\end{center}
subject to $t_*>0$ and $(a_0,a_1)\in\mathbb{R}^2$ satisfying (\ref{eq32}) at the end of Section \ref{sec5}.

We choose to present all our results in the one dimensional case, while many of these results can in fact be extended to convex Lipschitz or analytic functions $F(t_1,t_2,\cdots,t_m)$ of $m$ variables  naturally. Most of our results also hold on transitive subshift spaces of finite type, with some technical adjustments in their proofs involving the transition matrix.   

The organization of the work is as following. In Section \ref{sec2} we introduce some basics in thermodynamic formalism and the Central Limit Theorem for the process generated by a potential and the shift map on the symbolic space of finite type. We give an explicit bound on the tail term in the CLT. Section \ref{sec3} is devoted to the proof of Theorem \ref{thm8} and \ref{thm5}. We formulate some expression of the derivatives of the pressure (Corollary \ref{cor2}) linking directly to the CLT, which allows us to unveil the relationship between derivatives of the pressure function of various orders. Section \ref{sec4} is devoted to the proof of Theorem \ref{thm7}. In Section \ref{sec5} we consider fitting 1- and 2-level candidate analytic germs locally by pressure functions of locally constant potentials (Problem \ref{prob}) on symbolic spaces of finite type. We conjecture that any reasonable analytic germ of finite level can be fitted by the pressure function of some locally constant potential locally, as long as the number of the symbols is large enough.

\section{Thermodynamic formalism and the CLT}\label{sec2}

In this section we collect some basic notions and results in thermodynamic formalism for later use. We start from the pressure. Let $\Lambda$ be some symbolic set of finite symbols, $\Lambda^\mathbb{N}$ be the shift space equipped with the visual metric
\begin{center}
$d(x,y)=\cfrac{1}{2^{l(x,y)}}$
\end{center}
for distinct $x=x_0x_1x_2\cdots, y=y_0y_1y_2\cdots\in \Lambda^\mathbb{N}$, in which  
\begin{center}
$l(x,y)=\min\{i\in\mathbb{N}: x_i\neq y_i\}$. 
\end{center}
For a continuous potential $\phi: \Lambda^\mathbb{N}\rightarrow\mathbb{R}$ on the compact metric space $\Lambda^\mathbb{N}$, Let 
\begin{center}
$S_{m,\phi}(x)=\sum_{i=0}^{m-1}\phi\circ\sigma^i(x)$
\end{center}
for $m\in\mathbb{N}$, in which $\sigma$ is the shift map.
\begin{defn}\label{def1}
The pressure $P(\phi)$ of a continuous potential $\phi$ on $\Lambda^\mathbb{N}$ is defined to be
\begin{center}
$P(\phi)=\lim_{m\rightarrow\infty}\cfrac{1}{m}\log \sum_{\sigma^m(x)=x}e^{S_{m,\phi}(x)}$.
\end{center}
\end{defn}
One can refer to \cite[p208]{Wal1} for a definition for continuous potentials on general compact metric spaces. It satisfies the well-known variational formula
\begin{center}
$P(\phi)=\sup\{h(\mu)+\int\phi d\mu: \mu \mbox{ is a } \sigma-\mbox{invariant measure on }\Lambda^\mathbb{N}\}$.
\end{center}
Let $C^0(\Lambda^\mathbb{N})$ be the collection of all the continuous potentials on $\Lambda^\mathbb{N}$. Two potentials $\psi, \phi\in C^0(\Lambda^\mathbb{N})$ are said  to be \emph{cohomologous} \cite{Wal2} in case there exists a continuous map $\varphi: \Lambda^\mathbb{N}\rightarrow\mathbb{R}$ such that
\begin{center}
$\psi(x)-\phi(x)=\varphi(x)-\varphi\circ\sigma(x)$.
\end{center}
We write $\psi\sim\phi$ to denote the equivalence relationship between two potentials cohomologous to each other. The maps in 
\begin{center}
$\{\varphi(x)-\varphi\circ\sigma(x): \varphi\in C^0(\Lambda^\mathbb{N})\}$
\end{center}
are called \emph{coboundaries}. The importance of the cohomologous relationship is revealed in the following result.
\begin{pro}\label{pro1}
If $\psi\sim\phi$, then $P(\psi)=P(\phi)$. Moreover, $\psi$ and $\phi$ share the same equilibrium state.
\end{pro}
One can find a proof in \cite{Rue1} or \cite{PP}. Another important tool in thermodynamic formalism is the \emph{Ruelle operator}.

\begin{defn}
For a continuous potential $\psi: \Lambda^\mathbb{N}\rightarrow\mathbb{R}$, define the Ruelle operator $\mathcal{L}_\psi$ acting on $C^0(\Lambda^\mathbb{N})$ as
\begin{center}
$(\mathcal{L}_\psi f)(x)=\sum_{y:\sigma(y)=x}e^{\psi(y)}f(y)$
\end{center}
for $f\in C^0(\Lambda^\mathbb{N})$. 
\end{defn}

One can see easily that its compositions satisfy
\begin{equation}\label{eq13}
(\mathcal{L}^m_\psi f)(x)=\sum_{y:\sigma^m(y)=x}e^{S_{m,\psi}(y)}f(y)
\end{equation}
for any $m\in\mathbb{N}$. In case of $\psi$ being H\"older, it admits a simple maximum isolated eigenvalue $\lambda=e^{P(\psi)}$ such that,
\begin{equation}
(\mathcal{L}_\psi w_\psi)(x)=e^{P(\psi)}w_\psi(x)
\end{equation}  
for some eigenfunction $w_\psi(x)\in C^{0,h}(\Lambda^\mathbb{N})$, refer to \cite{Rue1}. The unique equilibrium measure for the H\"older potential $\psi$ is denoted by $\mu_\psi$. It then follows that
\begin{equation}\label{eq14}
(\mathcal{L}^m_\psi w_\psi)(x)=e^{mP(\psi)}w_\psi(x)
\end{equation}  
for $w_\psi(x)\in C^{0,h}(\Lambda^\mathbb{N})$. A potential $\psi$ is said to be \emph{normalized} if 
\begin{center}
$P(\psi)=0$ and $w_\psi=1_{\Lambda^\mathbb{N}}$,
\end{center}
in which $1_{\Lambda^\mathbb{N}}$ is the identity map on $\Lambda^\mathbb{N}$. In case of $\psi$ being not normalized, we call 
\begin{center}
$\bar{\psi}=\psi+\log w_\psi-\log w_\psi\circ\sigma-P(\psi)$
\end{center}
the \emph{normalization} of $\psi$. It is easy to check that $\bar{\psi}$ is a normalized potential. Moreover, $\bar{\psi}$ and $\psi$ share the same equilibrium state. 

Now we turn to the Central Limit Theorem for the random process $\{\phi\circ\sigma^j(x)\}_{j=0}^\infty$ with the equilibrium measure $\mu_\psi$ defined by some H\"older potential $\psi$, while $\phi$ is also assumed to be H\"older. It deals with the asymptotic behaviour of the distribution of $\cfrac{S_{m,\phi}}{\sqrt{m}}$ with respect to $\mu_\psi$ as $m\rightarrow\infty$. The Ruelle operator comes in here, see \cite{CP, Lal, Rou}. Let
\begin{center}
$G_m(y)=\mu_\psi\Big\{x\in \Lambda^\mathbb{N}: \cfrac{S_{m,\phi}(x)}{\sqrt{m}}<y\Big\}$
\end{center}
for $y\in\mathbb{R}$. For $a,b\in\mathbb{R}$ and $b>0$, Let $N_{a,b}(y)$ be the normal distribution with expectation $a$ and standard deviation $\sqrt{b}$ on $\mathbb{R}$, that is,
\begin{center}
$\cfrac{dN_{a,b}(y)}{dy}=\cfrac{1}{\sqrt{2\pi b}}e^{-(y-a)^2/2b}$
\end{center}
for $y\in\mathbb{R}$. For H\"older potentials $\psi, \phi$ on a shift space of finite type, since the pressure $P(\psi+t\phi)$ is analytic in a small neighbourhood around $0$, denote by
\begin{center}
$\Delta_m=P^{(m)}(\psi+t\phi)|_{t=0}$
\end{center}
for $m\in\mathbb{N}$ for convenience, while the readers can understand its dependence on $\psi, \phi$ easily from the contexts in the following. Let 
\begin{center}
$P(\psi+t\phi)=\sum_{m=0}^\infty\cfrac{\Delta_m}{m!}t^m=\sum_{m=0}^3\cfrac{\Delta_m}{m!}t^m+t^4\kappa(t)$,
\end{center}
in which $\kappa(t)=\sum_{m=0}^\infty\cfrac{\Delta_{m+4}}{(m+4)!}t^m$. 
\begin{CLT}
Let $\psi, \phi$ be H\"older potentials on a shift space of finite type with $\phi$ being not cohomologous to a constant. If $\int\phi d\mu_\psi=0$, we have
\begin{center}
$\lim_{m\rightarrow\infty} G_m(y)=N_{0,\Delta_2}(y)+O(1/\sqrt{m})$,
\end{center}
in which 
\begin{equation}\label{eq29}
O(1/\sqrt{m})\leq\cfrac{9|\Delta_3|+2|\Delta_4|}{\sqrt{2\pi^3 m}(\Delta_2)^{3/2}}.
\end{equation}
The convergence is uniform with respect to $y$.  In case of $\phi$ being generic, the result can be strengthened to
\begin{center}
$\lim_{m\rightarrow\infty} G_m(y)=N_{0,\Delta_2}(y)+H_m(y)+o(1/\sqrt{m})$,
\end{center}
in which $H_m(y)=\cfrac{\Delta_3}{6\sqrt{m}}\Big(1-\cfrac{y^2}{\Delta_2}\Big)e^{-\frac{y^2}{2\Delta_2}}$.
\end{CLT}

This fits into special cases of the \emph{Berry-Esseen Theorem} \cite{Fel}. There is nothing new in the version here comparing with \cite[Theorem 2, Theorem 3]{CP} or \cite[Theorem 4.13]{PP}, except the explicit bound on the tail term $O(1/\sqrt{m})$ in (\ref{eq29}). In the following we justify this explicit bound. To do this, let 
\begin{center}
$\chi_m(z)=\displaystyle\int e^{iz\frac{S_{m,\phi}}{\sqrt{m}}}d\mu_\psi$
\end{center}  
be the Fourier transformation of $G_m(y)$. Note that the Fourier transformation of $N_{0,\Delta_2}(y)$ is $e^{-\frac{z^2\Delta_2}{2}}$.
\begin{lemma}\label{lem1}
Let $\psi, \phi$ be H\"older potentials on a shift space of finite type with $\phi$ being not cohomologous to a constant. For $\epsilon>0$ small enough, we have
\begin{equation}\label{eq27}
\cfrac{1}{2\pi}\int_{0}^{\epsilon\sqrt{m}}\cfrac{1}{z}\Big|\chi_m(z)-e^{-\frac{z^2\Delta_2}{2}}\Big|dz\leq \cfrac{\sqrt{2} |\Delta_3|}{12\sqrt{\pi m}(\Delta_2)^{3/2}}
\end{equation} 
for any $m\in\mathbb{N}$ large enough.
\end{lemma}
\begin{proof}
According to \cite[(4.6)]{PP}, we have
\begin{center}
$\displaystyle\int_{0}^{\epsilon\sqrt{m}}\cfrac{1}{z}\Big|\chi_m(z)-e^{-\frac{z^2\Delta_2}{2}}+\cfrac{iz^3\Delta_3}{6\sqrt{m}}e^{-\frac{z^2\Delta_2}{2}}\Big|dz=O(1/m)$
\end{center}
for $\epsilon>0$ small enough. So
\begin{equation}\label{eq28}
\cfrac{1}{2\pi}\displaystyle\int_{0}^{\epsilon\sqrt{m}}\cfrac{1}{z}\Big|\chi_m(z)-e^{-\frac{z^2\Delta_2}{2}}\Big|dz\leq O(1/m)+\cfrac{|\Delta_3|}{12\pi\sqrt{m}}\displaystyle\int_{0}^{\epsilon\sqrt{m}}z^2e^{-\frac{z^2\Delta_2}{2}}dz.
\end{equation}
Considering
\begin{center}
$\displaystyle\int_{-\infty}^{\infty}z^2e^{-\frac{z^2\Delta_2}{2}}dz=\cfrac{\sqrt{2\pi}}{(\Delta_2)^{3/2}}$,
\end{center}
we obtain (\ref{eq27}) from (\ref{eq28}).
\end{proof}

Equipped with Lemma \ref{lem1} we can justify the explicit bound on the tail term in the  Central Limit Theorem in (\ref{eq29}).\\ 

Proof of the tail term in CLT:\\

\begin{proof}
Without loss of generality, suppose $\psi$ is normalized and $\int \phi d\mu_\psi=0$. It suffices for us to justify (\ref{eq29}) considering \cite[Theorem 2, Theorem 3]{CP}. Similar to the proof of \cite[Theorem 2]{CP}, apply \cite[Lemma 2]{Fel} with the cumulative functions $G_m(y)$ and $N_{0,\Delta_2}(y)$ in our case, one gets (c.f. \cite[(20)]{CP})
\begin{equation}\label{eq30}
|G_m(y)-N_{0,\Delta_2}(y)|\leq\cfrac{1}{2\pi}\int_{0}^{\epsilon\sqrt{m}}\cfrac{1}{z}\Big|\chi_m(z)-e^{-\frac{z^2\Delta_2}{2}}\Big|dz+\cfrac{24}{\epsilon\sqrt{2m\pi^3\Delta_2}}.
\end{equation}
Now let us take 
\begin{center}
$\cfrac{1}{\epsilon}=\cfrac{2}{\Delta_2}\big(\cfrac{|\Delta_3|}{6}+\cfrac{|\Delta_4|}{24}+\delta\big)$
\end{center}
for some small $\delta>0$, such that it satisfies (c.f. \cite[(10)]{CP}) 
\begin{center}
$\cfrac{1}{\epsilon}>\max\Big\{\cfrac{2}{\Delta_2}\Big(\cfrac{|\Delta_3|}{6}+t\kappa(t)\Big), \cfrac{2}{\Delta_2}\kappa(t)\Big\}$
\end{center}
for any $|t|<\epsilon$ in (\ref{eq30}). Considering (\ref{eq27}), we have

\begin{equation}\label{eq31}
\begin{array}{ll}
& |G_m(y)-N_{0,\Delta_2}(y)|\vspace{3mm}\\
\leq & \cfrac{\sqrt{2} |\Delta_3|}{12\sqrt{\pi m}(\Delta_2)^{3/2}}+\cfrac{24}{\sqrt{2m\pi^3\Delta_2}}\cfrac{2}{\Delta_2}\big(\cfrac{|\Delta_3|}{6}+\cfrac{|\Delta_4|}{24}+\delta\big)\vspace{3mm}\\
= & \cfrac{\sqrt{2} |\Delta_3|}{12\sqrt{\pi m}(\Delta_2)^{3/2}}+\cfrac{8|\Delta_3|}{\sqrt{2\pi^3 m}(\Delta_2)^{3/2}}+\cfrac{2|\Delta_4|}{\sqrt{2\pi^3 m}(\Delta_2)^{3/2}}+\cfrac{48\delta}{\sqrt{2\pi^3 m}(\Delta_2)^{3/2}}\vspace{3mm}\\
\leq & \cfrac{9|\Delta_3|}{\sqrt{2\pi^3 m}(\Delta_2)^{3/2}}+\cfrac{2|\Delta_4|}{\sqrt{2\pi^3 m}(\Delta_2)^{3/2}}+\cfrac{48\delta}{\sqrt{2\pi^3 m}(\Delta_2)^{3/2}}.
\end{array}
\end{equation}
Finally, let $\delta\rightarrow 0$ in (\ref{eq31}), we get (\ref{eq29}).

\end{proof}

We will deal with the pressure function $P(\psi+t\phi)$ for $t\geq 0$ and $\psi, \phi\in C^{0,h}(\Lambda^\mathbb{N})$ for some $0<h\leq 1$ in the following sections. By \cite{Rue2}, $P(\psi+t\phi)$ depends analytically on $t$ in case that $\psi, \phi$ are H\"older. We will often assume that
\begin{center}
$\int\phi d\mu_\psi=0$
\end{center}
in the following when dealing with the higher differentials of $P(\psi+t\phi)$ because if $\int\phi d\mu_\psi=c\neq 0$, we have 
\begin{center}
$P(\psi+t(\phi-c))=P(\psi+t\phi)-ct$, 
\end{center}
then
\begin{equation}\label{eq18}
\cfrac{d^n P(\psi+t(\phi-c))}{dt^n}=\cfrac{d^n P(\psi+t\phi)}{dt^n}
\end{equation}
for any $n\geq 2$ while  $\int(\phi-c)d\mu_\psi=0$. We can also assume that $\psi$ is normalized when dealing with the differentials of $P(\psi+t\phi)$. If this is not the case we can simply change $\psi$ to its normalization $\bar{\psi}$ while
\begin{equation}\label{eq21}
\cfrac{d^n P(\psi+t\phi)}{dt^n}=\cfrac{d^n P(\bar{\psi}+t\phi)}{dt^n}
\end{equation}
for $n\geq 1$ because 
\begin{center}
$P(\bar{\psi}+t\phi)=P(\psi+t\phi)-P(\psi)$
\end{center}
for any $t\in\mathbb{R}$.

\section{Derivatives of the pressures of H\"older potentials}\label{sec3}

In this section we formulate some explicit expressions for the derivatives of the pressure $P(t\phi)=P(t)$ in terms of the derivatives of the eigenfunction of $\mathcal{L}_{t\phi}$ for $\phi\in C^{0,h}(\Lambda^\mathbb{N})$ with respect to $t$. We give basically two expressions of the derivatives, one of which allows the introduction of the random stochastic process $\{\phi\circ\sigma^j(x)\}_{j=0}^m$ for $m\in\mathbb{N}$. Upon the expression  we prove Theorem \ref{thm8} and \ref{thm5} in virtue of the CLT for the random process $\{\phi\circ\sigma^j(x)\}_{j=0}^\infty$.

First we define some basics to deal with the higher derivatives of compositional functions by the \emph{Fa\`a di Bruno's formula}. For an integer $j\in\mathbb{N}$, we say 
\begin{center}
$\tau=\tau_1\tau_2\cdots\tau_q$
\end{center}
with $q\in\mathbb{N}$ is a \emph{partition} of $j$ if the non-increasing sequence of positive integers $j\geq\tau_1\geq\tau_2\geq\cdots\geq\tau_q\geq 1$ satisfies $\sum_{i=1}^q\tau_i=j$. Denote the collection of all the possible partitions of $j$ by $\mathfrak{P}(j)$. For example, Table \ref{tab1} lists all the partitions in $\mathfrak{P}(5)$.

\begin{table}[ht]
\caption{Partitions of $5$} 
\centering 
\begin{tabular}{c c} 
\hline 
5 & q=1\\
4,1 & q=2\\
3,2 & q=2\\
3,1,1 & q=3\\
2,2,1 &q=3\\ 
2,1,1,1 &q=4\\
1,1,1,1,1 &q=5\\
\hline 
\end{tabular}
\label{tab1} 
\end{table}

We sometimes simply write $\tau$ to denote the set $\{\tau_1,\tau_2,\cdots,\tau_q\}$ for convenience in the following, so $\#\tau=q$. Now for $\tau$ being a partition of $j\geq 1$, let $\{B_j^\tau\}$ be the number of different choices of dividing a set of $j$ different elements into $\#\tau=q$ sets of sizes $\{\tau_i\}_{i=1}^q$ respectively (with no order on the sets of partitions). Set $B_0^0=1$ for convenience. For example, consider the cases $j=5$ and $\tau=3,1,1$, the number of different choices of dividing a set of $5$ different elements into $q=3$ sets of sizes $3,1,1$ respectively is 
\begin{center}
$C_5^3=10=B_5^{3,1,1}$. 
\end{center}
Table \ref{tab2} lists all the numbers $\{B_5^\tau\}_{\tau\in\mathfrak{P}(5)}$. 

\begin{table}[ht]
\caption{The coefficients $B_5^\tau$} 
\centering 
\begin{tabular}{c c} 
\hline 
$B_5^5=1$\\
$B_5^{4,1}=5$\\
$B_5^{3,2}=10$\\
$B_5^{3,1,1}=10$\\
$B_5^{2,2,1}=15$\\ 
$B_5^{2,1,1,1}=10$\\
$B_5^{1,1,1,1,1}=1$\\
\hline 
\end{tabular}
\label{tab2} 
\end{table}

For a smooth map $f: X\rightarrow Y$ between two metric spaces $X, Y$ and some partition $\tau=\tau_1,\tau_2,\cdots,\tau_q\in\mathfrak{P}(j)$ with $j\geq 1$, let
\begin{center}
$f^{(\tau)}(x)=f^{(\tau_1)}(x)f^{(\tau_2)}(x)\cdots f^{(\tau_q)}(x)$
\end{center}
be the product of the derivatives.  For $j=0$ and $\tau=0\in\mathfrak{P}(0)$, set $f^{(0)}(x)=1$.  Then for two smooth functions $f: X\rightarrow Y$ and $g:Y\rightarrow Z$ between metric spaces $X, Y, Z$, we have
\begin{equation}\label{eq5}
\cfrac{d^j(g\circ f(x))}{dx^j}=\sum_{\tau\in\mathfrak{P}(j)}B_j^\tau g^{(\#\tau)}(f(x))f^{(\tau)}(x) 
\end{equation}    
in virtue of Fa\`a di Bruno's formula.

Now we turn to the higher differentials of the pressure function.  We start by considering some standard case, then extend the result to the general case.

\begin{theorem}\label{thm2}
Let $\psi, \phi\in C^{0,h}(\Lambda^\mathbb{N})$ with $\psi$ being normalized for some finite symbolic set $\Lambda$. Assume $\int \phi d\mu_{\psi}=0$, in which $\mu_{\psi}$ is the equilibrium state of $\psi$. Let $w(t,x)$ be the eigenfunction of the maximum isolated eigenvalue $e^{P(\psi+t\phi)}$ of $\mathcal{L}_{\psi+t\phi}$, which depends analytically on $t$ in a small neighbourhood of $0$. Considering the differentials of the pressure function $P(\psi+t\phi)$ at $t=0$, we have
\begin{equation}\label{eq8}
\begin{array}{ll}
P^{(n)}(\psi+t\phi)|_{t=0}=&\sum_{j=1}^n C_n^j\int_{\Lambda^{\mathbb{N}}}(\phi(x))^j w^{(n-j)}(0,x)d\mu_{\psi}(x)\\
&-\sum_{j=2}^{n-2}C_n^j\sum_{\tau\in\mathfrak{P}(j),1\notin\tau} B_j^\tau P^{(\tau)}(\psi+t\phi)|_{t=0}\int_{\Lambda^{\mathbb{N}}}w^{(n-j)}(0,x)d\mu_{\psi}(x)\\
&-\sum_{\tau\in\mathfrak{P}(n),\{1,n\}\cap\tau=\emptyset} B_n^\tau P^{(\tau)}(\psi+t\phi)|_{t=0}
\end{array}
\end{equation}
for any $n\geq 2$.
\end{theorem}

\begin{proof}
According to the above notations, note that
\begin{equation}
(\mathcal{L}_{\psi+t\phi} w(t,\cdot))(x)=e^{P(\psi+t\phi)}w(t,x).
\end{equation}  

The $n$-th derivative of $(\mathcal{L}_{\psi+t\phi} w(t,\cdot))(x)=\sum_{y:\sigma(y)=x}e^{\psi(y)+t\phi(y)}w(t,y)$ gives
\begin{equation}\label{eq2}
\begin{array}{ll}
&\cfrac{d^n \mathcal{L}_{\psi+t\phi} w(t,\cdot)(x)}{dt^n}\vspace{1mm}\\
=&\sum_{y:\sigma(y)=x}\sum_{j=0}^{n}C_n^j \cfrac{d^j e^{(\psi+t\phi)(y)}}{dt^j}w^{(n-j)}(t,y)\vspace{3mm}\\
=&\sum_{y:\sigma(y)=x}\sum_{j=0}^{n}C_n^j e^{(\psi+t\phi)(y)}(\phi(y))^j w^{(n-j)}(t,y)\vspace{3mm}\\
=&\sum_{j=0}^{n}C_n^j \mathcal{L}_{\psi+t\phi}\big((\phi(\cdot))^j w^{(n-j)}(t,\cdot)\big).
\end{array}
\end{equation}
All differentials are with respect to $t$. In case of $t=0$ this means
\begin{equation}\label{eq4} 
\cfrac{d^n \mathcal{L}_{\psi+t\phi} w(t,\cdot)(x)}{dt^n}|_{t=0}=\sum_{j=0}^{n}C_n^j \mathcal{L}_{\psi}\big((\phi(\cdot))^j w^{(n-j)}(0,\cdot)\big).
\end{equation}
Note that the dual operator $\mathcal{L}_{\psi}^*$ fixes $\mu_{\psi}$,  so integration of both sides of (\ref{eq4}) gives
\begin{equation}\label{eq7}
\int\cfrac{d^n \mathcal{L}_{\psi+t\phi} w(t,\cdot)(x)}{dt^n}|_{t=0}d\mu_{\psi}(x)=\sum_{j=0}^{n}C_n^j \int(\phi(x))^j w^{(n-j)}(0,x)\mu_{\psi}(x).
\end{equation} 

In order to get the $n$-th derivative of $P(\psi+t\phi)$, differentiate $e^{P(\psi+t\phi)}w(t,x)$ for $n$ times by (\ref{eq5}), we get
\begin{equation}\label{eq3}
\begin{array}{ll}
&\cfrac{d^n \Big(e^{P(\psi+t\phi)}w(t,x)\Big)}{dt^n}\vspace{1mm}\\
=&\sum_{j=0}^{n}C_n^j \cfrac{d^j e^{P(\psi+t\phi)}}{dt^j}w^{(n-j)}(t,x)\vspace{3mm}\\
=&\sum_{j=0}^{n-1}C_n^j \cfrac{d^j e^{P(\psi+t\phi)}}{dt^j}w^{(n-j)}(t,x)+\cfrac{d^n e^{P(\psi+t\phi)}}{dt^n}w(t,x)\vspace{3mm}\\
=&\sum_{j=0}^{n-1}C_n^j \sum_{\tau\in\mathfrak{P}(j)} B_j^\tau P^{(\tau)}(\psi+t\phi)e^{P(\psi+t\phi)}w^{(n-j)}(t,x)\vspace{3mm}\\
&+\sum_{\tau\in\mathfrak{P}(n)} B_n^\tau P^{(\tau)}(\psi+t\phi)e^{P(\psi+t\phi)}w(t,x)\vspace{3mm}\\
=&\sum_{j=0}^{n-1}C_n^j \Big(\sum_{\tau\in\mathfrak{P}(j),1\notin\tau} B_j^\tau P^{(\tau)}(\psi+t\phi)+\sum_{\tau\in\mathfrak{P}(j),1\in\tau} B_j^\tau P^{(\tau)}(\psi+t\phi)\Big)e^{P(\psi+t\phi)}w^{(n-j)}(t,x)\vspace{3mm}\\
&+\sum_{\tau\in\mathfrak{P}(n),n\notin\tau} B_n^\tau P^{(\tau)}(\psi+t\phi)e^{P(\psi+t\phi)}w(t,x)+P^{(n)}(\psi+t\phi)e^{P(\psi+t\phi)}w(t,x).\vspace{3mm}\\
\end{array}
\end{equation} 

Remember $P(\psi)=0$ and $w(0,x)=1$ as $\psi$ is normalized (\cite[p66]{PP}). Take $t=0$ in (\ref{eq3}) we get
\begin{equation}\label{eq49}
\begin{array}{ll}
&\cfrac{d^n \Big(e^{P(\psi+t\phi)}w(t,x)\Big)}{dt^n}|_{t=0}\vspace{1mm}\\
=&\sum_{j=0}^{n-1}C_n^j \Big(\sum_{\tau\in\mathfrak{P}(j),1\notin\tau} B_j^\tau P^{(\tau)}(\psi+t\phi)|_{t=0}+\sum_{\tau\in\mathfrak{P}(j),1\in\tau} B_j^\tau P^{(\tau)}(\psi+t\phi)|_{t=0}\Big)w^{(n-j)}(0,x)\vspace{3mm}\\
&+\sum_{\tau\in\mathfrak{P}(n),n\notin\tau} B_n^\tau P^{(\tau)}(\psi+t\phi)|_{t=0}+P^{(n)}(\psi+t\phi)|_{t=0}\vspace{3mm}\\
\end{array}
\end{equation}

Since $\int \phi d\mu_{\psi}=P'(\psi+t\phi)|_{t=0}=0$ and $\int w'(0,x)d\mu_{\psi}=0$ (\cite[p66]{PP}), integrate both sides of (\ref{eq49}) with respect to $\mu_{\psi}$, we get 
\begin{equation}\label{eq6}
\begin{array}{ll}
&\displaystyle\int\cfrac{d^n \Big(e^{P(\psi+t\phi)}w(t,x)\Big)}{dt^n}|_{t=0}d\mu_{\psi}\vspace{1mm}\\
=&\sum_{j=0}^{n-1}C_n^j \sum_{\tau\in\mathfrak{P}(j),1\notin\tau} B_j^\tau P^{(\tau)}(\psi+t\phi)|_{t=0}\int w^{(n-j)}(0,x)d\mu_{\psi}\vspace{3mm}\\
&+\sum_{\tau\in\mathfrak{P}(n),\{1,n\}\cap\tau=\emptyset} B_n^\tau P^{(\tau)}(\psi+t\phi)|_{t=0}+P^{(n)}(\psi+t\phi)|_{t=0}.\vspace{3mm}\\
\end{array}
\end{equation}

Finally, combining (\ref{eq7}) and (\ref{eq6}) together we get (\ref{eq8}).

\end{proof}

\begin{rem}
The terms 
\begin{center}
$-\sum_{j=2}^{n-2}C_n^j\sum_{\tau\in\mathfrak{P}(j),1\notin\tau} B_j^\tau P^{(\tau)}(\psi+t\phi)|_{t=0}\int_{\Lambda^{\mathbb{N}}}w^{(n-j)}(0,x)d\mu_{\psi}(x)$
\end{center}
and
\begin{center}
$-\sum_{\tau\in\mathfrak{P}(n),\{1,n\}\cap\tau=\emptyset} B_n^\tau P^{(\tau)}(\psi+t\phi)|_{t=0}$
\end{center}
in (\ref{eq8}) are null in case of $n\leq 3$. This also applies to the corresponding terms later.
\end{rem}

\begin{rem}\label{rem1}
These appear to be inductive formulas, while one can always get non-inductive ones via substituting the lower differentials $P^{(\tau)}(\psi+t\phi)|_{t=0}$ by their non-inductive versions depending only on $\phi(x), \{w^{(j)}(0,x)\}_{j=1}^n$ and $\mu_{\psi}(x)$. This also applies to Theorem \ref{thm4}.
\end{rem}

One can find some description of derivatives of the pressure function by \emph{covariance} of the sequence of functions $\{\phi\circ\sigma^j\}_{j\in\mathbb{N}}$ in \cite[Corollary 1]{KS1} for smooth $\phi$.  Without the assumptions of $\psi$ being normalized and $\int \phi d\mu_{\psi}=0$,  Theorem \ref{thm2} evolves into the following form.

\begin{coro}\label{cor1}
Let $\psi, \phi\in C^{0,h}(\Lambda^\mathbb{N})$ with  some finite symbolic set $\Lambda$. $\mathcal{L}_{\psi+t\phi}$ admits a maximum isolated eigenvalue $e^{P(\psi+t\phi)}$ close to $e^{P(\psi)}$ with eigenfunction $w(t,x)$ whose projection depends analytically on $t$ in a small neighbourhood of $0$. Considering the differentials of the pressure $P(\psi+t\phi)$ at $t=0$, we have
\begin{equation}\label{eq10}
\begin{array}{ll}
P^{(n)}(\psi+t\phi)|_{t=0}=&\sum_{j=1}^n C_n^j\displaystyle\int_{\Lambda^{\mathbb{N}}}\Big(\phi(x)-\int \phi d\mu_{\psi}\Big)^j w^{(n-j)}(0,x)d\mu_{\psi}(x)\vspace{3mm}\\
&-\sum_{j=2}^{n-2}C_n^j\sum_{\tau\in\mathfrak{P}(j),1\notin\tau} B_j^\tau P^{(\tau)}(\psi+t\phi)|_{t=0}\displaystyle\int_{\Lambda^{\mathbb{N}}}w^{(n-j)}(0,x)d\mu_{\psi}(x)\vspace{3mm}\\
&-\sum_{\tau\in\mathfrak{P}(n),\{1,n\}\cap\tau=\emptyset} B_n^\tau P^{(\tau)}(\psi+t\phi)|_{t=0}
\end{array}
\end{equation}
for any $n\geq 2$.
\end{coro}

\begin{proof}
Let
\begin{center}
$\bar{\psi}=\psi+\log w_\psi(x)-\log w_\psi\circ\sigma-P(\psi)$
\end{center}
in which $w_\psi(x)$ is the eigenfunction of $\mathcal{L}_{\psi}$ corresponding to the eigenvalue $e^{P(\psi)}$.  Take pressure in the following equation
\begin{center}
$\bar{\psi}+t\phi=\psi+t\phi+\log w_\psi(x)-\log w_\psi\circ\sigma-P(\psi)$,
\end{center}
then apply Proposition \ref{pro1}, we see that
\begin{center}
$P(\bar{\psi}+t\phi)=P(\psi+t\phi)-P(\psi)$.
\end{center} 
This implies
\begin{equation}\label{eq9}
\cfrac{d^n P(\bar{\psi}+t\phi)}{dt^n}=\cfrac{d^nP(\psi+t\phi)}{dt^n}
\end{equation}
for any $n\geq 1$. Now apply Theorem \ref{thm2} to the normalized potential $\bar{\psi}$ and $\phi-\int \phi d\mu_{\psi}$, (note that $\int \big(\phi-\int \phi d\mu_{\psi}\big)d\mu_{\psi}=0$ and $\mu_{\psi}=\mu_{\bar{\psi}}$), we justify the corollary considering (\ref{eq9}).
\end{proof}

In the following we present some concrete formulas of some special order $n$ in virtue of Theorem \ref{thm2} for later use. 
\begin{coro}
Let $\psi, \phi\in C^{0,h}(\Lambda^\mathbb{N})$ with $\psi$ being normalized. Let $\mu_{\psi}$ be the equilibrium state of $\psi$ and $\int \phi d\mu_{\psi}=0$. Let $e^{P(\psi+t\phi)}$ be the maximum eigenvalue of $\mathcal{L}_{\psi+t\phi}$ with eigenfunction $w(t,x)$ for small $t$. Then we have
\begin{equation}\label{eq11}
P^{'''}(\psi+t\phi)|_{t=0}=3\int \phi w^{''}(0,x)d\mu_{\psi}+3\int \phi^2 w^{'}(0,x)d\mu_{\psi}+\int \phi^3 d\mu_{\psi}.
\end{equation}
\end{coro}
\begin{proof}
This follows instantly from Theorem \ref{thm2} with $n=3$, along with some direct computations on the Fa\`a di Bruno's coefficients $\{B_3^\tau\}_{\tau\in\mathfrak{P}(3)}$.
\end{proof}

\begin{coro}
Let $\psi, \phi\in C^{0,h}(\Lambda^\mathbb{N})$ with $\psi$ being normalized. Let $\mu_{\psi}$ be the equilibrium state of $\psi$ and $\int \phi d\mu_{\psi}=0$. Let $e^{P(\psi+t\phi)}$ be the maximum eigenvalue of $\mathcal{L}_{\psi+t\phi}$ with eigenfunction $w(t,x)$ for small $t$. Then we have
\begin{equation}\label{eq12}
\begin{array}{ll}
& P^{''''}(\psi+t\phi)|_{t=0}\vspace{3mm}\\
=& 4\int \phi w^{'''}(0,x)d\mu_{\psi}+6\int \phi^2 w^{''}(0,x)d\mu_{\psi}+4\int \phi^3 w^{'}(0,x)d\mu_{\psi}+\int \phi^4 d\mu_{\psi}\vspace{3mm}\\
&-6P^{''}(\psi+t\phi)|_{t=0}\int w^{''}(0,x)d\mu_{\psi}-3(P^{''}(\psi+t\phi)|_{t=0})^2\vspace{3mm}\\
=& 4\int \phi w^{'''}(0,x)d\mu_{\psi}+6\int \phi^2 w^{''}(0,x)d\mu_{\psi}+4\int \phi^3 w^{'}(0,x)d\mu_{\psi}+\int \phi^4 d\mu_{\psi}\vspace{3mm}\\
&-6(\int \phi^2d\mu_{\psi}+2\int \phi w'(0,x)d\mu_{\psi})\int w^{''}(0,x)d\mu_{\psi}-3(\int \phi^2d\mu_{\psi}+2\int \phi w'(0,x)d\mu_{\psi})^2.
\end{array}
\end{equation}

\end{coro}
\begin{proof}
The first equality follows instantly from Theorem \ref{thm2} with $n=4$ along with some direct computations on the Fa\`a di Bruno's coefficients $\{B_4^\tau\}_{\tau\in\mathfrak{P}(4)}$. The second one is true as
\begin{center}
$P^{''}(\psi+t\phi)|_{t=0}=\int \phi^2d\mu_{\psi}+2\int \phi w'(0,x)d\mu_{\psi}.$
\end{center}
The latter description depends only on $\phi(x), \{w^{(j)}(0,x)\}_{j=1}^3$ and $\mu_{\psi}(x)$.
\end{proof}

One can also get some precise formulas for some particular $n$ in Corollary \ref{cor1}, and some non-inductive ones as we indicate in Remark \ref{rem1}. While the formulas (\ref{eq8}, \ref{eq10}, \ref{eq11}, \ref{eq12}) all give interesting descriptions of the differentials of the pressure function $P(\psi+t\phi)$, it seems to us difficult to discover any essential rigid restriction on them, or relationships between them.  In the following we turn to the description of them by the random stochastic process $\{\phi\circ\sigma^j(x)\}_{j=0}^\infty$. This is not a new idea on exploring the regularity of the pressure function $P(\psi+t\phi)$, as one can recall it from many others' work in thermodynamic formalism. Again we first consider some standard case, then extend to the general case. 

\begin{theorem}\label{thm4}
Let $\psi, \phi\in C^{0,h}(\Lambda^\mathbb{N})$ with $\psi$ being normalized. Let $\mu_{\psi}$ be the equilibrium state of $\psi$ and $\int \phi d\mu_{\psi}=0$. Let $e^{P(\psi+t\phi)}$ be the maximum isolated eigenvalue of $e^{P(\psi+t\phi)}$ with eigenfunction $w(t,x)$ whose projection depends analytically on $t$. Considering the differentials of the pressure $P(\psi+t\phi)$ at $t=0$, we have
\begin{equation}\label{eq15}
\begin{array}{ll}
&P^{(n)}(\psi+t\phi)|_{t=0}\\
=&\lim_{m\rightarrow\infty}\cfrac{1}{m}\Big(\sum_{j=2}^n C_n^j\displaystyle\int_{\Lambda^\mathbb{N}}(S_{m,\phi}(x))^j w^{(n-j)}(0,x)d\mu_{\psi}(x)\vspace{3mm}\\
&-\sum_{j=2}^{n-2}C_n^j\sum_{\tau\in\mathfrak{P}(j),1\notin\tau} m^{\#\tau}B_j^\tau P^{(\tau)}(\psi+t\phi)|_{t=0}\displaystyle\int_{\Lambda^\mathbb{N}}w^{(n-j)}(0,x)d\mu_{\psi}(x)\vspace{3mm}\\
&-\sum_{\tau\in\mathfrak{P}(n),\{1,n\}\cap\tau=\emptyset} m^{\#\tau} B_n^\tau P^{(\tau)}(\psi+t\phi)|_{t=0}\Big)
\end{array}
\end{equation}
for any $n\geq 2$.
\end{theorem}

\begin{proof}
The proof follows the routine of Proof of Theorem \ref{thm2}. Considering  (\ref{eq13}), we do $n$-differentials on both sides of (\ref{eq14}), take $t=0$, then integrate both sides with respect to $\mu_{\psi}(x)$, divided by $m$, we get 
\begin{equation}\label{eq16}
\begin{array}{ll}
&P^{(n)}(\psi+t\phi)|_{t=0}\\
=&\cfrac{1}{m}\Big(\sum_{j=1}^n C_n^j\displaystyle\int_{\Lambda^\mathbb{N}}(S_{m,\phi}(x))^j w^{(n-j)}(0,x)d\mu_{\psi}(x)\vspace{3mm}\\
&-\sum_{j=2}^{n-2}C_n^j\sum_{\tau\in\mathfrak{P}(j),1\notin\tau} m^{\#\tau}B_j^\tau P^{(\tau)}(\psi+t\phi)|_{t=0}\displaystyle\int_{\Lambda^\mathbb{N}}w^{(n-j)}(0,x)d\mu_{\psi}(x)\vspace{3mm}\\
&-\sum_{\tau\in\mathfrak{P}(n),\{1,n\}\cap\tau=\emptyset} m^{\#\tau} B_n^\tau P^{(\tau)}(\psi+t\phi)|_{t=0}\Big)
\end{array}
\end{equation}  
as (\ref{eq8}). Now since $w^{(n-1)}(0,x)$ is bounded on $X$, the ergodic theorem guarantees
\begin{equation}\label{eq17}
\lim_{m\rightarrow\infty}\cfrac{1}{m}\int_{\Lambda^\mathbb{N}}S_{m,\phi}(x) w^{(n-1)}(0,x)d\mu_{\psi}(x)=0.
\end{equation}
Then (\ref{eq15}) follows from (\ref{eq16}) as $m\rightarrow\infty$ considering (\ref{eq17}).
\end{proof}

Theorem \ref{thm4} establishes some link between the differentials of the pressure function and the the process $\{\phi\circ\sigma^j(x)\}_{j=0}^\infty$ through $S_{m,\phi}$ with respect to the equilibrium state $\mu_{\psi}$. We also formulate a general version of the result.
  
\begin{coro}
Let $\psi, \phi\in C^{0,h}(\Lambda^\mathbb{N})$ with  $\mu_{\psi}$ be the equilibrium state of $\psi$. $\mathcal{L}_{\psi+t\phi}$ admits a maximum isolated eigenvalue $e^{P(\psi+t\phi)}$ close to $e^{P(\psi)}$ with eigenfunction $w(t,x)$ whose projection depends analytically on $t$ in a small neighbourhood of $0$. Considering the differentials of the pressure function $P(\psi+t\phi)$ at $t=0$, we have
\begin{equation}
\begin{array}{ll}
&P^{(n)}(\psi+t\phi)|_{t=0}\\
=&\lim_{m\rightarrow\infty}\cfrac{1}{m}\Big(\sum_{j=2}^n C_n^j\displaystyle\int_{\Lambda^\mathbb{N}}\big(S_{m,\phi}-m\int\phi d\mu_{\psi}\big)^j w^{(n-j)}(0,x)d\mu_{\psi}(x)\vspace{3mm}\\
&-\sum_{j=2}^{n-2}C_n^j\sum_{\tau\in\mathfrak{P}(j),1\notin\tau} m^{\#\tau}B_j^\tau P^{(\tau)}(\psi+t\phi)|_{t=0}\displaystyle\int_{\Lambda^\mathbb{N}}w^{(n-j)}(0,x)d\mu_{\psi}(x)\vspace{3mm}\\
&-\sum_{\tau\in\mathfrak{P}(n),\{1,n\}\cap\tau=\emptyset} m^{\#\tau} B_n^\tau P^{(\tau)}(\psi+t\phi)|_{t=0}\Big)
\end{array}
\end{equation}
for any $n\geq 2$.
\end{coro}
\begin{proof}
Equipped with Theorem \ref{thm4}, the proof follows in line with the Proof of Corollary \ref{cor1}.
\end{proof}

Now we give some precise descriptions of the third and fourth differentials of $P(\psi+t\phi)$ in virtue of Theorem \ref{thm4}.
\begin{coro}
Let $\psi, \phi\in C^{0,h}(\Lambda^\mathbb{N})$ with $\psi$ being normalized. Let $\mu_{\psi}$ be the equilibrium state of $\psi$ and $\int \phi d\mu_{\psi}=0$. Let $e^{P(\psi+t\phi)}$ be the maximum eigenvalue of $\mathcal{L}_{\psi+t\phi}$ with eigenfunction $w(t,x)$ for small $t$. Then we have
\begin{equation}\label{eq25}
P^{''}(\psi+t\phi)|_{t=0}=\lim_{m\rightarrow\infty}\cfrac{1}{m}\int S_{m,\phi}^2d\mu_{\psi}.
\end{equation}
\end{coro}

\begin{rem}
$P^{''}(\psi+t\phi)|_{t=0}$ is called variance of the random process $\{\phi\circ\sigma^j(x)\}_{j=0}^\infty$, whose name can be interpreted from the Central Limit Theorem. See \cite{Rue1, PP}.
\end{rem}

\begin{coro}\label{cor2}
Let $\psi, \phi\in C^{0,h}(\Lambda^\mathbb{N})$ with $\psi$ being normalized. Let $\mu_{\psi}$ be the equilibrium state of $\psi$ and $\int \phi d\mu_{\psi}=0$. Let $e^{P(\psi+t\phi)}$ be the maximum eigenvalue of $\mathcal{L}_{\psi+t\phi}$ with eigenfunction $w(t,x)$ for small $t$. Then we have
\begin{equation}\label{eq22}
P^{'''}(\psi+t\phi)|_{t=0}=\lim_{m\rightarrow\infty}\cfrac{3}{m}\int S_{m,\phi}^2 w^{'}(0,x)d\mu_{\psi}+\lim_{m\rightarrow\infty}\cfrac{1}{m}\int S_{m,\phi}^3 d\mu_{\psi}.
\end{equation}
\end{coro}
\begin{proof}
This follows instantly from Theorem \ref{thm4} with $n=3$.
\end{proof}

\begin{coro}\label{cor5}
Let $\psi, \phi\in C^{0,h}(\Lambda^\mathbb{N})$ with $\psi$ being normalized. Let $\mu_{\psi}$ be the equilibrium state of $\psi$ and $\int \phi d\mu_{\psi}=0$. Let $e^{P(\psi+t\phi)}$ be the maximum eigenvalue of $\mathcal{L}_{\psi+t\phi}$ with eigenfunction $w(t,x)$ for small $t$. Then we have
\begin{equation}
\begin{array}{ll}
& P^{(4)}(\psi+t\phi)|_{t=0}\vspace{3mm}\\
=& \lim_{m\rightarrow\infty}\Big(\cfrac{6}{m}\displaystyle\int S_{m,\phi}^2 w^{''}(0,x)d\mu_{\psi}+\cfrac{4}{m}\int S_{m,\phi}^3 w^{'}(0,x)d\mu_{\psi}+\cfrac{1}{m}\int S_{m,\phi}^4 d\mu_{\psi}\vspace{3mm}\\
&-6P^{''}(\psi+t\phi)|_{t=0}\displaystyle\int w^{''}(0,x)d\mu_{\psi}-3m(P^{''}(\psi+t\phi)|_{t=0})^2\Big)\vspace{3mm}\\
=& \lim_{m\rightarrow\infty}\Big(\cfrac{6}{m}\displaystyle\int S_{m,\phi}^2 w^{''}(0,x)d\mu_{\psi}+\cfrac{4}{m}\int S_{m,\phi}^3 w^{'}(0,x)d\mu_{\psi}+\cfrac{1}{m}\int S_{m,\phi}^4 d\mu_{\psi}\vspace{3mm}\\
&-\cfrac{6}{m}\displaystyle\int S_{m,\phi}^2 d\mu_{\psi}\displaystyle\int w^{''}(0,x)d\mu_{\psi}-\cfrac{3}{m}(\int S_{m,\phi}^2 d\mu_{\psi})^2\Big)\vspace{3mm}\\.
\end{array}
\end{equation}

\end{coro}
\begin{proof}
The first equality follows instantly from Theorem \ref{thm4} with $n=4$, while the second one is true considering (\ref{eq25}). The last description depends only on $\phi(x), \{w^{(j)}(0,x)\}_{j=1}^2$ and $\mu_{\psi}(x)$.
\end{proof}

Through the above formulas we see the importance of the asymptotic distribution of the random variable $S_{m,\phi}$ with respect to $\mu_{\psi}$, which is describe by the Central Limit Theorem for the process $\{\phi\circ\sigma^j(x)\}_{j=0}^\infty$. Equipped with all the above results, now we are in a position to prove the rigidity results on the third differentials of $P(\psi+t\phi)$ upon Corollary \ref{cor2}. We first show Theorem \ref{thm5}.

\begin{proof}[Proof of Theorem \ref{thm5}]
From now on we fix $t_*\in (0,\infty)$. Let $\psi=t_*\phi$.  Simply by making a change of variable we can see that
\begin{center}
$P^{(n)}(t_*)=P^{(n)}(t\phi)|_{t=t_*}=P^{(n)}(\psi+t\phi)|_{t=0}$
\end{center} 
for any $n\geq 0$. So (\ref{eq19}) is equivalent to 
\begin{equation}\label{eq20}
|P'''(\psi+t\phi)|_{t=0}\big(1-\sqrt{2\pi}(P''(\psi+t\phi)|_{t=0})^{3/2}\big)|\leq 3 M_\phi P''(\psi+t\phi)|_{t=0}.
\end{equation}
We can assume $\psi$ is normalized as otherwise we can change it to its normalization considering (\ref{eq21}). Moreover, it suffices for us to prove it under the assumption $\int\phi d\mu_\psi=0$ in virtue of (\ref{eq18}). If $P''(\psi+t\phi)|_{t=0}=0$, then $\phi$ is cohomologous to a constant according to \cite[Proposition 4.12]{PP}. This forces $P'''(\psi+t\phi)|_{t=0}=0$, so (\ref{eq20}) is satisfied in this case. In the following we assume $P''(\psi+t\phi)|_{t=0}>0$. We resort to Corollary \ref{cor2} to justify (\ref{eq20}) under the above assumptions. We first estimate the term $\cfrac{1}{m}\int S_{m,\phi}^3 d\mu_{\psi}$ in (\ref{eq22}). Now the Central Limit Theorem comes in. 
\begin{center}
$
\begin{array}{ll}
&\cfrac{1}{m}\displaystyle\int S_{m,\phi}^3 d\mu_{\psi}\vspace{3mm}\\
=& \sqrt{m}\displaystyle\int \Big(\cfrac{S_{m,\phi}}{\sqrt{m}}\Big)^3 d\mu_{\psi}\vspace{3mm}\\
=& \sqrt{m}\displaystyle\int y^3 d G_m(y)\vspace{3mm}\\
=& \sqrt{m}\displaystyle\int y^3 d N_{0,P''(\psi+t\phi)|_{t=0}}(y)+\sqrt{m}\int y^3 dH_m(y)+\sqrt{m}\cdot o(1/\sqrt{m})\vspace{3mm}\\
=& \sqrt{m}\cdot 0+\displaystyle\int y^3 d\Big(\cfrac{P'''(\psi+t\phi)|_{t=0}}{6}\big(1-\cfrac{y^2}{P''(\psi+t\phi)|_{t=0}}\big)e^{-y^2/2P''(\psi+t\phi)|_{t=0}}\Big)\\
&+\sqrt{m}\cdot o(1/\sqrt{m})\vspace{3mm}\\
=& P'''(\psi+t\phi)|_{t=0}\sqrt{2\pi}(P''(\psi+t\phi)|_{t=0})^{3/2}+\sqrt{m}\cdot o(1/\sqrt{m}).\vspace{3mm}\\
\end{array}
$
\end{center}
By taking $m\rightarrow\infty$ we get
\begin{equation}
\lim_{m\rightarrow\infty} \cfrac{1}{m}\displaystyle\int S_{m,\phi}^3 d\mu_{\psi}=P'''(\psi+t\phi)|_{t=0}\sqrt{2\pi}(P''(\psi+t\phi)|_{t=0})^{3/2}.
\end{equation}
Considering (\ref{eq22}) we have
\begin{equation}\label{eq23}
P'''(\psi+t\phi)|_{t=0}\big(1-\sqrt{2\pi}(P''(\psi+t\phi)|_{t=0})^{3/2}\big)=\lim_{m\rightarrow\infty}\cfrac{3}{m}\int S_{m,\phi}^2 w^{'}(0,x)d\mu_{\psi}.
\end{equation}
Since $w^{'}(0,x)$ depends continuously on $x\in X$, there exists some $M_\phi$ depending on $\phi$, such that 
\begin{equation}\label{eq24}
|w^{'}(0,x)|\leq M_\phi. 
\end{equation}
Now taking absolute values on both sides of (\ref{eq23}) we justify (\ref{eq20}), considering (\ref{eq24}) and (\ref{eq25}).

\end{proof}

The proof of Theorem \ref{thm8} on the pressure functions of non-generic H\"older potentials follows a similar way.

\begin{proof}[Proof of Theorem \ref{thm8}]
Fix $t_*\in (0,\infty)$, we can simply assume $\psi=t_*\phi$ is normalised and $\int\phi d\mu_\psi=0$. In case that $P''(\psi+t\phi)|_{t=0}=0$, so $\phi$ is cohomologous to a constant, (\ref{eq50}) holds obviously. In the following we assume $\phi$ is not cohomologous to a constant, so $P''(\psi+t\phi)|_{t=0}>0$. We again resort to Corollary \ref{cor2} to justify (\ref{eq50}) under these assumptions. Now for the term $\cfrac{1}{m}\int S_{m,\phi}^3 d\mu_{\psi}$ in (\ref{eq22}), in virtue of the Central Limit Theorem, 
\begin{equation}\label{eq51}
\begin{array}{ll}
&\cfrac{1}{m}\displaystyle\int S_{m,\phi}^3 d\mu_{\psi}\vspace{3mm}\\
=& \sqrt{m}\displaystyle\int \Big(\cfrac{S_{m,\phi}}{\sqrt{m}}\Big)^3 d\mu_{\psi}\vspace{3mm}\\
=& \sqrt{m}\displaystyle\int y^3 d G_m(y)\vspace{3mm}\\
\leq& \sqrt{m}\displaystyle\int y^3 d N_{0,P''(\psi+t\phi)|_{t=0}}(y)+\sqrt{m}\cfrac{9|P'''(\psi+t\phi)|_{t=0}|+2|P^{(4)}(\psi+t\phi)|_{t=0}|}{\sqrt{2\pi^3 m}(P''(\psi+t\phi)|_{t=0})^{3/2}}\vspace{3mm}\\
=& \sqrt{m}\cdot 0+\cfrac{9|P'''(\psi+t\phi)|_{t=0}|+2|P^{(4)}(\psi+t\phi)|_{t=0}|}{\sqrt{2\pi^3}(P''(\psi+t\phi)|_{t=0})^{3/2}}\vspace{3mm}\\
=& \cfrac{9|P'''(\psi+t\phi)|_{t=0}|+2|P^{(4)}(\psi+t\phi)|_{t=0}|}{\sqrt{2\pi^3}(P''(\psi+t\phi)|_{t=0})^{3/2}}\\
\end{array}
\end{equation}
for $m$ large enough. By taking $m\rightarrow\infty$ in (\ref{eq51}), we get
\begin{equation}\label{eq52}
\lim_{m\rightarrow\infty} \cfrac{1}{m}\displaystyle\int S_{m,\phi}^3 d\mu_{\psi}\leq \cfrac{9|P'''(\psi+t\phi)|_{t=0}|+2|P^{(4)}(\psi+t\phi)|_{t=0}|}{\sqrt{2\pi^3}(P''(\psi+t\phi)|_{t=0})^{3/2}}.
\end{equation}
Taking modulus on both sides of (\ref{eq52}) we get
\begin{equation}
\begin{array}{ll}
& |P'''(\psi+t\phi)|_{t=0}|\vspace{3mm}\\
\leq & \Big|\lim_{m\rightarrow\infty}\cfrac{1}{m}\displaystyle\int S_{m,\phi}^3 d\mu_{\psi}\Big|+\Big|\lim_{m\rightarrow\infty}\cfrac{3}{m}\int S_{m,\phi}^2 w^{'}(0,x)d\mu_{\psi}\Big|\vspace{3mm}\\
\leq & \cfrac{9|P'''(\psi+t\phi)|_{t=0}|+2|P^{(4)}(\psi+t\phi)|_{t=0}|}{\sqrt{2\pi^3}(P''(\psi+t\phi)|_{t=0})^{3/2}}+3M_\phi \Big|\lim_{m\rightarrow\infty}\cfrac{3}{m}\displaystyle\int S_{m,\phi}^2 d\mu_{\psi}\Big|\vspace{3mm}\\
= & \cfrac{9|P'''(\psi+t\phi)|_{t=0}|+2|P^{(4)}(\psi+t\phi)|_{t=0}|}{\sqrt{2\pi^3}(P''(\psi+t\phi)|_{t=0})^{3/2}}+3M_\phi P''(\psi+t\phi)|_{t=0}
\end{array}
\end{equation}
for some $|w^{'}(0,x)|\leq M_\phi$, which results in (\ref{eq50}).
\end{proof}

One can predict from Corollary \ref{cor5}, Theorem \ref{thm4} and the proof of Theorem \ref{thm8}, Theorem \ref{thm5} that some more rigid relationships between higher differentials of the pressure function $\{P^{(n)}(t\phi)\}_{n\in\mathbb{N}}$ are possible. These rigidity relationships impose restrictions on fitting convex analytic functions whose supporting lines intersecting the vertical axis in some bounded set in $[0,\infty)$ by pressures of H\"older potentials.

\section{Global Fitting of convex analytic functions via pressures of H\"older potentials}\label{sec4}

This section is dedicated to the proof of Theorem \ref{thm7}. We start from the following result on some global behaviour of the pressure functions of generic H\"older potentials.

\begin{theorem}\label{thm6}
Let $\alpha>0$. If a strictly convex analytic function $F(t)$ on $(\alpha,\infty)$, with its supporting lines intersecting the vertical axis in $[\underline{\gamma},\overline{\gamma}]\subset [0,\infty)$, such that
\begin{equation}\label{eq26}
\sup_{t\in(\alpha,\infty)}\Bigg\{\Big|\cfrac{F^{'''}(t)}{F^{''}(t)}-\sqrt{2\pi F^{''}(t)}\Big|\Bigg\}=\infty,
\end{equation}
then there does not exist any generic H\"older potential $\phi$ on any shift space $X$ of finite type satisfying
\begin{center}
$P(t\phi)= F(t)$
\end{center}
on $(\alpha,\infty)$.
\end{theorem}

\begin{proof}
This follows directly from Theorem \ref{thm5} in fact. Suppose on the contrary that there exist some shift space $X$ of finite type and some generic H\"older potential $\phi\in C^{0,h}(X)$ satisfying $P(t\phi)= F(t)$ on $(\alpha,\infty)$, then according to Theorem \ref{thm5}, we have
\begin{center}
$\sup_{t\in(\alpha,\infty)}\Bigg\{\Big|\cfrac{F^{'''}(t)}{F^{''}(t)}-\sqrt{2\pi F^{''}(t)}
\Big|\Bigg\}\leq 3M_\phi$
\end{center}
for some finite $M_\phi>0$. This contradicts (\ref{eq26}). 
\end{proof}

Be careful that we cannot exclude the possibility that one can locally fit some convex analytic function through the pressure of some generic H\"older potential on some shift space of finite type by Theorem \ref{thm5}. This is because for any strictly convex analytic function $F(t)$ on $(\alpha,\infty)$ and $\alpha\leq\underline{\alpha}\leq\overline{\alpha}$, we always have 
\begin{center}
$\sup_{\underline{\alpha}\leq t\leq \overline{\alpha}}\Bigg\{\Big|\cfrac{F^{'''}(t)}{F^{''}(t)}-\sqrt{2\pi F^{''}(t)}
\Big|\Bigg\}<\infty$.
\end{center}
So one cannot exclude the possibility that there exists some generic H\"older potential $\phi$ on some shift space of finite type satisfying
\begin{center}
$P(t\phi)= F(t)$
\end{center}
on $[\underline{\alpha},\overline{\alpha}]$ through Theorem \ref{thm5}. See Section \ref{sec5} for more results on the problem of local fitting of some convex analytic functions through the pressures of H\"older potentials.

Now for $\alpha>0$, let
\begin{center}
$
\begin{array}{ll}
\mathcal{F}_\alpha=&\{F(t): F(t) \mbox{ is a strictly convex analytic function on } (\alpha,\infty) \mbox{ satisfying } (\ref{eq26}), \vspace{3mm}\\
& \mbox{ its supporting lines intersect the vertical axis in a bounded interval in } [0,\infty)\}.\vspace{3mm}
\end{array}
$
\end{center}

We will show that $\mathcal{F}_\alpha\neq\emptyset$ for any $\alpha>0$ in the following. 

\begin{pro}\label{pro3}
For any $\alpha>0$, we have
\begin{center}
$\tilde{\mathcal{F}}_\alpha=\Big\{F_{a,b,c}(t)=\cfrac{at^2+bt+te^{-ct^2}+e^{-ct^2}}{t}\Big|_{(\alpha,\infty)}\Big\}_{a,b>0, c>1/2\sqrt{2}}\subset \mathcal{F}_\alpha$.
\end{center} 
\end{pro}
\begin{proof}
The restricted functions on $(\alpha,\infty)$ are of course analytic. Considering the second derivative of a function $F_{a,b,c}(t)\in\tilde{\mathcal{F}}_\alpha$, we have
\begin{center}
$F_{a,b,c}^{''}(t)=4c^2t^2e^{-ct^2}+4c^2te^{-ct^2}-2ce^{-ct^2}+2ct^{-1}e^{-ct^2}+2t^{-3}e^{-ct^2}$
\end{center} 
for $t\in(0,\infty)$. Now since
\begin{center}
$4c^2t+2ct^{-1}\geq 2\sqrt{8c^3}>2c$
\end{center}
considering $c>1/2\sqrt{2}$, we can see that $F_{a,b,c}^{''}(t)>0$ on $(0,\infty)$. This shows that for any $\alpha>0$, $F_{a,b,c}(t)\in\tilde{\mathcal{F}}_\alpha$ is a convex function. Considering the third differential of a function $F_{a,b,c}(t)\in\tilde{\mathcal{F}}_\alpha$, we have
\begin{center}
$F_{a,b,c}^{'''}(t)=-8c^3t^3e^{-ct^2}-8c^3t^2e^{-ct^2}+12c^2te^{-ct^2}-6ct^{-2}e^{-ct^2}-6t^{-4}e^{-ct^2}$
\end{center} 
for $t\in(0,\infty)$. Then we have
\begin{center}
$\lim_{t\rightarrow\infty} \Bigg(\cfrac{F^{'''}_{a,b,c}(t)}{F^{''}_{a,b,c}(t)}-\sqrt{2\pi F^{''}_{a,b,c}(t)}\Bigg)=\lim_{t\rightarrow\infty} \cfrac{-8c^3t^3e^{-ct^2}}{4c^2t^2e^{-ct^2}}=-\infty$.
\end{center}
This means that $F_{a,b,c}(t)\in\tilde{\mathcal{F}}_\alpha$ satisfies (\ref{eq26}). To see that the supporting lines of  a function $F_{a,b,c}(t)\in\tilde{\mathcal{F}}_\alpha$ intersect the vertical axis in a bounded domain in $[0,\infty)$, write the function as 
\begin{center}
$F_{a,b,c}(t)=at+b+e^{-ct^2}+t^{-1}e^{-ct^2}$.
\end{center}
Its graph on $(0,\infty)$ is a strictly convex smooth curve with asymptotes $y=at+b$ and $t=0$.
\end{proof}

\begin{figure}[h]
\centering
\includegraphics[scale=1]{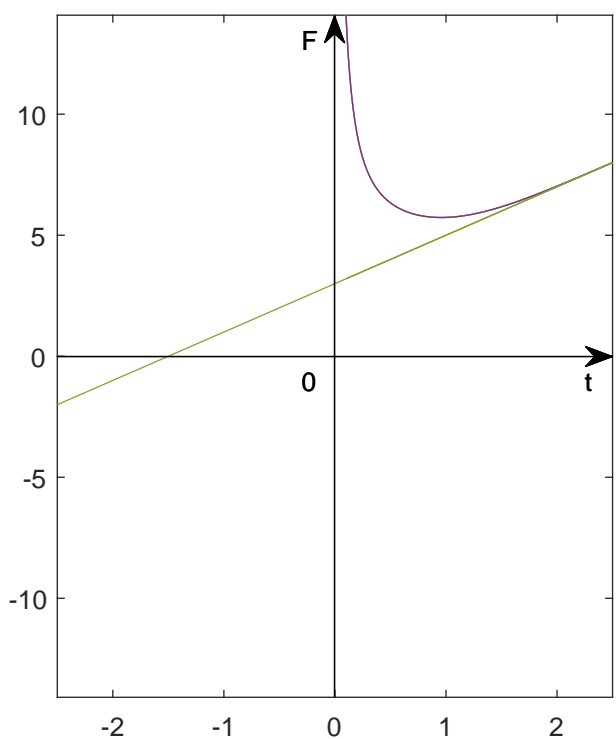}
\caption{ Graph of $F_{2,3,1}(t)$}
\label{fig1}
\end{figure}

In Figure \ref{fig1} we provide the readers with the graph of the function 
\begin{center}
$F_{2,3,1}(t)=\cfrac{2t^2+3t+te^{-t^2}+e^{-t^2}}{t}$
\end{center} 
on $(0,\infty)$.

This means that any function in the family $\tilde{\mathcal{F}}_\alpha$ cannot be fitted by any generic H\"older potential on any shift space of finite type globally, considering Theorem \ref{thm6}. In the following we deny the possibility that they can be fitted by \emph{non-generic} H\"older potentials on shift spaces of finite type. 

\begin{defn}
A continuous potential $\phi: X \to \mathbb R$ on a shift space $X$ of finite type is said to be non-generic if for some normalised potential $\psi$, the spectral radius of the complex Ruelle operator $\mathcal{L}_{\psi+it\phi}$ equals $1$ for some $t\neq 0$.
\end{defn}

One can show that if $\phi$ is non-generic then there exists a continuous function $u: X \to \mathbb R$, $c_\phi \in \mathbb R$ and a locally constant potential $\tilde{\phi}: X \to \mathbb R$, such that
\begin{equation}\label{eq54}
\phi(x) = u\circ\sigma(x)-u(x) + c_\phi + \tilde{\phi}(x).
\end{equation}

\begin{pro}\label{pro4}
For any $\alpha>0$ and any $F_{a,b,c}(t)\in\tilde{\mathcal{F}}_\alpha$ with $a,b>0, c>\frac{1}{2\sqrt{2}}$, there does not exist any non-generic H\"older potential $\phi$ on any shift space of finite type such that
\begin{center}
$P(t\phi)= F(t)$
\end{center}
on $(\alpha,\infty)$.
\end{pro}
\begin{proof}
Note that for a non-generic H\"older potential $\phi$ on a shift space of finite type, according to (\ref{eq54}), we have
\begin{center}
$P(t\phi)=tc_\phi  + P(t\tilde{\phi})$,
\end{center} 
in which $\tilde{\phi}$ is some locally constant potential. By the explicit formula (see Lemma \ref{lem2}) for the pressure functions of locally constant potentials on shift spaces of finite type, we see that any $F_{a,b,c}(t)$ cannot be fitted by pressure of any non-generic H\"older potential $\phi$ globally.
\end{proof}

Equipped with all the above results, Theorem \ref{thm7} follows instantly from Proposition \ref{pro3} and \ref{pro4}.

\section{Local fitting of analytic germs via pressures of locally constant potentials}\label{sec5}

In this section we deal with the local fitting of analytic functions by the pressures of H\"older potentials, especially the pressures of piecewise constant ones. Firstly we borrow some notion originating from analytic continuation.

\begin{defn}
A germ at $t_*$ is the sum of infinite power series
\begin{center}
$g(t)=a_0+a_1(t-t_*)+\cfrac{a_2}{2!}(t-t_*)^2+\cfrac{a_3}{3!}(t-t_*)^3+\cdots$
\end{center}
for some $(a_0,a_1,\cdots)\in\mathbb{R}^{\infty}$. 
\end{defn}

The convergent radius (the superior of values $\delta\geq 0$ on $[t_*-\delta,t_*+\delta]$ the germ converges) of the power series is called the \emph{radius} of the germ. We are only interested in germs of radius $\delta>0$. The following problem will be our concern in this section.

\begin{prob}\label{prob}
For a germ 
\begin{center}
$g(t)=a_0+a_1(t-t_*)+\cfrac{a_2}{2!}(t-t_*)^2+\cdots$ 
\end{center}
at $t_*$ with some strictly positive radius, does there exist some H\"older potential $\phi$ on some shift space of finite type and some $\delta>0$, such that
\begin{center}
$P(t\phi)=g(t)$
\end{center}
on $[t_*-\delta,t_*+\delta]$?
\end{prob}

The question can still be understood in Katok's flexibility program in the class of symbolic dynamical systems, or even in some smooth systems. Obvious conditions to guarantee a positive answer to the problem are (\ref{eq32}) and
\begin{equation}\label{eq33}
a_2>0.
\end{equation}
Condition (\ref{eq33}) guarantees convexity of the germ (in some neighbourhood of $t_*$) while (\ref{eq32}) guarantees the supporting lines of the germ intersect the vertical axis in a bounded set in $[0,\infty)$ (also in some neighbourhood of $t_*$). We are especially interested in its answer when the H\"older potential in Problem \ref{prob} is required to be a piecewise constant one. We have seen the importance of the family of locally constant potentials in approximating convex analytic functions in Corollary \ref{thm1}. In fact Corollary \ref{thm1} has some interesting interpretation in approximation theory \cite{Tim}, when we consider the explicit expressions of the pressures of locally constant potentials on the shift space of finite type. For $n\in\mathbb{N}$, recall that 
\begin{center}
$\Lambda_n=\{1,2,\cdots,n\}$.
\end{center} 
\begin{lemma}\label{lem2}  
For an integer $k\geq 0$, consider some locally constant potential
\begin{center}
$\phi(x)=c_{x_{-k}x_{-k+1}\cdots x_0\cdots x_{k-1}x_{k}}$
\end{center}
for $x=\cdots x_{-1}x_{0}x_{1}\cdots\in[x_{-k}\cdots x_{k}]$ on the shift space $\Lambda_n^{\mathbb{Z}}$, we have
\begin{center}
$P(t\phi)=\log\sum_{(x_{-k},\cdots,x_{k})\in\Lambda_n^{2k+1}}e^{tc_{x_{-k}\cdots x_{k}}}$
\end{center}
for any $t\in(-\infty,\infty)$.
\end{lemma}

\begin{proof}
This follows from \cite[Theorem 9.6]{Wal1} by some direct calculations through Definition \ref{def1} of the pressure.  See also \cite[p214]{Wal1}.
\end{proof}

\begin{rem}
The result can be extended to transitive subshifts of finite type. In this case the pressure is the logarithm of the maximal eigenvalue of some appropriate matrix.
\end{rem}

Now combining Corollary \ref{thm1} and Lemma  \ref{lem2}, we have the following result.

\begin{coro}\label{cor3}
Let $F(t)$ be a convex Lipschitz function on $(\alpha,\infty)$ for some $\alpha>0$, such that its supporting lines intersect the vertical axis in $[\underline{\gamma},\overline{\gamma}]$ with $0\leq\underline{\gamma}\leq\overline{\gamma}<\infty$. Then there exists some $K\in\mathbb{N}$ and some sequences of constants
\begin{center}
$\{c_{n,j}\}_{j=1}^{K^n}$,
\end{center}
such that
\begin{equation}\label{eq34}
\lim_{n\rightarrow\infty} \log \sum_{j=1}^{K^n} e^{tc_{n,j}}=F(t)
\end{equation}
for any $t\in(\alpha,\infty)$.
\end{coro}
\begin{proof}
Take $K=\#\Lambda$ for the symbolic set in the proof of Corollary \ref{thm1}, then the locally constant potential $\phi_{n}(x)=\phi_{n,-}(x)$ admits $K^n$ constant values respectively on corresponding level-$n$ cylinder sets. Denote these values by $\{c_{n,j}\}_{j=1}^{K^n}$ for $n\in\mathbb{N}$. According to Lemma \ref{lem2},
\begin{center}
$P(t\phi_{n,-})=\log \sum_{j=1}^{K^n} e^{tc_{n,j}}$
\end{center}
for any  $n\geq 1$. This gives (\ref{eq34}) by virtue of  (\ref{eq1}).

\end{proof}

Corollary \ref{cor3} indicates that logarithm of the finite sums of the exponential maps in the family $\{e^{tc}\}_{c\in\mathbb{R}}$ are dense in the space of certain convex Lipschitz maps on $(\alpha,\infty)$.  The above approximation is uniform with respect to $t$ in a bounded set. This makes the family $\{e^{tc}\}_{c\in\mathbb{R}}$ (family of locally constant potentials) important in detecting the properties of certain convex Lipschitz maps (among continuous or H\"older potentials). 

From now on we turn our attention to Problem \ref{prob}, but with restriction on locally constant potentials. We focus on locally constant potentials defined on the level-$0$ cylinder sets, whose theory is presumably parallel to the ones defined on the deeper cylinder sets. On the shift space $\Lambda_n^{\mathbb{Z}}$ with $n\geq 2$, consider the locally constant potential 
\begin{center}
$\phi(x)=z_{x_0}$
\end{center}
for $x=\cdots x_{-1}x_{0}x_{1}\cdots\in[x_0]$, in which $\{z_i\}_{1\leq i\leq n}$ are all constants. Let 
\begin{center}
$Q_0(t,z_1,z_2,\cdots,z_n)=\sum_{i=1}^n e^{tz_i}$,
\end{center}
so 
\begin{center}
$P(t\phi)=\log Q_0(t,z_1,\cdots,z_n)$ 
\end{center}
by Lemma \ref{lem2}. Let
\begin{center}
$Q_1(t,z_1,z_2,\cdots,z_n)=\sum_{i=1}^n z_ie^{tz_i}$
\end{center} 
and
\begin{center}
$Q_2(t,z_1,z_2,\cdots,z_n)=\sum_{1\leq i<j\leq n}(z_i-z_j)^2e^{t(z_i+z_j)}$.
\end{center} 
Through some elementary calculations one can check that
\begin{center}
$P'(t\phi)=\cfrac{dP(t\phi)}{dt}=\cfrac{Q_1(t,z_1,\cdots,z_n)}{Q_0(t,z_1,\cdots,z_n)}$
\end{center} 
while
\begin{equation}\label{eq39}
P''(t\phi)=\cfrac{d^2P(t\phi)}{dt^2}=\cfrac{Q_2(t,z_1,\cdots,z_n)}{Q_0^2(t,z_1,\cdots,z_n)}.
\end{equation}
Let
\begin{center}
$R_2(t,z_1,z_2,\cdots,z_n)=\sum_{i=1}^n z_i^2e^{tz_i}$,
\end{center} 
one can check that
\begin{center}
$Q_2(t,z_1,\cdots,z_n)=Q_0(t,z_1,\cdots,z_n)R_2(t,z_1,\cdots,z_n)-Q_1^2(t,z_1,\cdots,z_n)$.
\end{center}
In the following we will often fix $t=t_*>0$, so we will frequently write
\begin{center}
$Q_0(t_*,z_1,z_2,\cdots,z_n)=Q_0(z_1,z_2,\cdots,z_n)$
\end{center}
with $t_*$ omitted for convenience. Similar notations apply to other terms above. Let
\begin{equation}\label{eq35}
Q_0(z_1,\cdots,z_n)=\sum_{i=1}^n e^{t_*z_i}=e^{a_0},
\end{equation}
\begin{equation}\label{eq36}
Q_1(z_1,\cdots,z_n)=\sum_{i=1}^n z_ie^{t_*z_i}=a_1e^{a_0}
\end{equation}
be two equations with unknowns $\{z_1, z_2, \cdots, z_n\}$ for fixed $t_*>0, (a_0, a_1)\in\mathbb{R}^2$ and some $n\geq 2$. Let 
\begin{center}
$\Gamma_{\ref{eq35}}^n=\{(z_1, z_2, \cdots, z_n)\in\mathbb{R}^n: z_1, z_2, \cdots, z_n \mbox{ satisfy } (\ref{eq35})\}$
\end{center}
and
\begin{center}
$\Gamma_{\ref{eq36}}^n=\{(z_1, z_2, \cdots, z_n)\in\mathbb{R}^n: z_1, z_2, \cdots, z_n \mbox{ satisfy } (\ref{eq36})\}$.
\end{center}

They are both $n-1$ dimensional smooth hypersurfaces.  We first present readers with the following result on fitting an analytic function
\begin{center}
$a_0+a_1(t-t_*)+O((t-t_*)^2)$
\end{center} 
with $t_*, a_0, a_1$ subject to (\ref{eq32}) around some fixed $t_*>0$ by pressures of locally constant potentials on general shift spaces of finite type.

\begin{theorem}\label{thm9}
Let $t_*>0, (a_0, a_1)\in\mathbb{R}^2, n\geq 2$ satisfying (\ref{eq32}) and 
\begin{equation}\label{eq38}
\cfrac{a_0-\log n}{t_*}<a_1.
\end{equation} 
Then there exists some $\delta_n>0$ and some sequence $\{r_{i,n}\}_{i=1}^n\subset \mathbb{R}$, such that the locally constant potential 
\begin{center}
$\phi(x)=r_{x_0}$
\end{center}
for $x=\cdots x_{-1}x_{0}x_{1}\cdots\in[x_0]$ on the full shift space $\Lambda_n^{\mathbb{Z}}$ satisfies
\begin{center}
$P(t\phi)=a_0+a_1(t-t_*)+O((t-t_*)^2)$
\end{center}
on $[t_*-\delta_n,t_*+\delta_n]$.
\end{theorem}
\begin{proof}
In fact it suffices for us to show that the system of equations
\begin{center}
$
\left\{
\begin{array}{ll}
(\ref{eq35}),\\
(\ref{eq36})
\end{array}
\right.$
\end{center}
with unknowns $\{z_1, z_2, \cdots, z_n\}$ admits a solution under conditions of the theorem. Without loss of generality we assume 
\begin{equation}\label{eq37}
z_1\leq z_2\leq \cdots\leq z_n.
\end{equation} 
Under this assumption, it is easy to see that
\begin{center}
$\cfrac{a_0-\log n}{t_*}\leq z_n <\cfrac{a_0}{t_*}$.
\end{center}
Now we estimate the values of $Q_1(z_1,\cdots,z_n)$ with $z_n$ approaching the terminals. When $z_n$ approaches the right terminal from below, we have
\begin{center}
$\lim_{(z_1, z_2, \cdots, z_n)\in\Gamma_{\ref{eq35}}^n,\ z_n\nearrow\frac{a_0}{t_*}} Q_1(z_1,\cdots,z_n)=\cfrac{a_0}{t_*}e^{a_0}>a_1e^{a_0}$
\end{center}
in virtue of (\ref{eq32}). When $z_n$ approaches the left terminal from above, we have
\begin{center}
$\lim_{(z_1, z_2, \cdots, z_n)\in\Gamma_{\ref{eq35}}^n,\ z_n\searrow\frac{a_0}{t_*}} Q_1(z_1,\cdots,z_n)=\cfrac{a_0-\log n}{t_*}e^{a_0}<a_1e^{a_0}$
\end{center}
in virtue of (\ref{eq38}). Since $\Gamma_{\ref{eq35}}^n$ is a smooth hypersurface, by the mean value theorem, there exists some $(r_{1,n}, r_{2,n}, \cdots, r_{n,n})\in\Gamma_{\ref{eq35}}^n$ satisfying (\ref{eq35}) and (\ref{eq36}) simultaneously. At last, for $x=\cdots x_{-1}x_{0}x_{1}\cdots\in[x_0]$ on the full shift space $\Lambda_n^{\mathbb{Z}}$, let  
\begin{center}
$\phi(x)=r_{x_0,n}$
\end{center}
be the locally constant potential. As $P(t\phi)$ is analytic, there exists some $\delta_n>0$ such that 
\begin{center}
$P(t\phi)=a_0+a_1(t-t_*)+O((t-t_*)^2)$
\end{center}
for $t\in[t_*-\delta_n,t_*+\delta_n]$.
\end{proof}

\begin{rem}
The core step in the proof of Theorem \ref{thm9} is in fact finding the extremes of the function $Q_1(z_1,\cdots,z_n)$ subject to  (\ref{eq35}), (\ref{eq32}) and (\ref{eq38}). One can detect the points of extremes by the Karush-Kuhn-Tucker (KKT) conditions \cite{Kar, KT}, which generalizes the method of Lagrange multipliers by allowing inequality subjections. 
\end{rem}

Be careful that those $\{r_{i,n}\}_{i=1}^n$ all depend on $n$ in fact.  Theorem \ref{thm9} induces the following interesting flexibility result on fitting certain analytic functions locally by pressures of locally constant potentials on general shift space of finite type.

\begin{coro}\label{cor4}
Let $t_*>0$ and $(a_0, a_1)\in\mathbb{R}^2$ satisfy (\ref{eq32}). Then there exists some $N\in\mathbb{N}$, such that for any $n\geq N$, there exist some some $\delta_n>0$ and some sequence $\{r_{i,n}\}_{i=1}^n\subset \mathbb{R}$, such that the locally constant potential 
\begin{center}
$\phi(x)=r_{x_0,n}$
\end{center}
for $x=\cdots x_{-1}x_{0}x_{1}\cdots\in[x_0]$ on the full shift space $\Lambda_n^{\mathbb{Z}}$ satisfies
\begin{center}
$P(t\phi)=a_0+a_1(t-t_*)+O((t-t_*)^2)$
\end{center}
on $[t_*-\delta_n,t_*+\delta_n]$. 
\end{coro}
\begin{proof}
Under conditions of the corollary, for the given values $t_*, a_0, a_1$ satisfying (\ref{eq32}), choose $N\in\mathbb{N}$ large enough such that 
\begin{center}
$\cfrac{a_0-\log N}{t_*}<a_1$.
\end{center}
This means that for any $n>N$ condition (\ref{eq38}) is satisfied for $t_*, a_0, a_1, n$. Then the conclusion follows from Theorem \ref{thm9}.
\end{proof}

Note that on some particular symbolic spaces Theorem \ref{thm9} and \ref{cor4} may be trivial. For example, for given $(t_*, a_0, a_1)\in\mathbb{R}^3$ without any subjections, by choosing $\beta=e^{a_0-t_*a_1}$, consider the constant potential
\begin{center}
$\phi(x)= a_1$
\end{center}
on the $\beta$-shift space with symbols $\{0,1,\cdots, \lfloor\beta\rfloor\}$. It is easy to see that 
\begin{center}
$P(t\phi)=a_0-t_*a_1+a_1 t=a_0+a_1(t-t_*)$
\end{center}
on $(-\infty,\infty)$. However, our results guarantee conclusions on general shift spaces.

From now on we go towards the proof of Theorem \ref{thm3}.  For fixed $t_*>0, (a_0, a_1)\in\mathbb{R}^2$ and $n\geq 3$, let 
\begin{center}
$\Gamma_{\ref{eq35},\ref{eq36}}^n=\Gamma_{\ref{eq35}}^n\cap\Gamma_{\ref{eq36}}^n=\{(z_1, z_2, \cdots, z_n)\in\mathbb{R}^n: z_1, z_2, \cdots, z_n \mbox{ satisfy } (\ref{eq35}) \mbox{ and } (\ref{eq36})\}$.
\end{center}
We describe some topological properties of the set $\Gamma_{\ref{eq35},\ref{eq36}}^n$ in the following result.

\begin{figure}[h]
\centering
\includegraphics[scale=0.95]{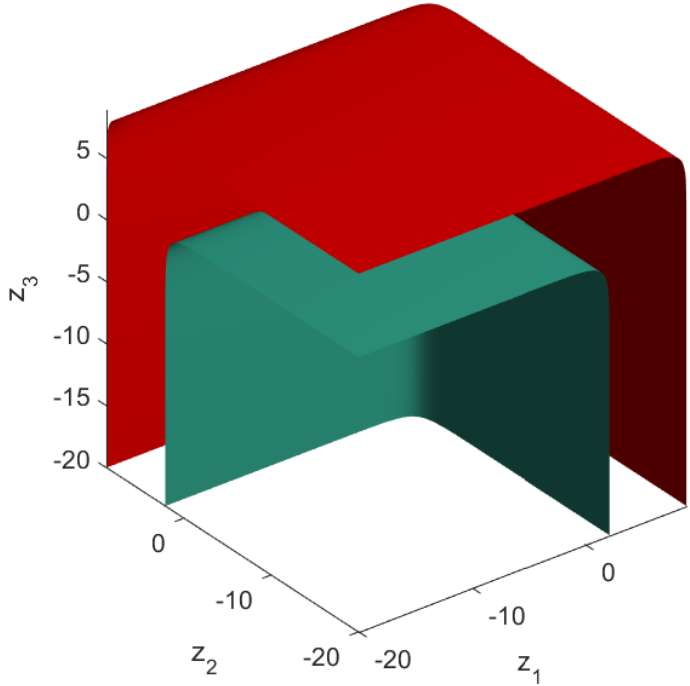}
\caption{$\Gamma_{\ref{eq35},1,2,1}^3$ (green) and $\Gamma_{\ref{eq36},1,2,1}^3$ (red)}
\label{fig2}
\end{figure}

\begin{lemma}
For fixed $t_*>0, (a_0, a_1)\in\mathbb{R}^2$  subject to (\ref{eq32}) and $n\geq 3$, in case $\Gamma_{\ref{eq35},\ref{eq36}}^n\neq\emptyset$ and $a_1\neq \cfrac{a_0-\log n}{t_*}$, it is a compact $(n-2)$-dimension smooth manifold. 
\end{lemma}
\begin{proof}
The Jacobian of the functions $Q_0(z_1,\cdots,z_n)-e^{a_0}$ and $Q_1(z_1,\cdots,z_n)-a_1e^{a_0}$ with respect to $z_1, z_2, \cdots, z_n$ is
\begin{center}
$
J=\begin{pmatrix} t_*e^{t_*z_1} & t_*e^{t_*z_2} & \cdots & t_*e^{t_*z_n}\\ 
e^{t_*z_1}+t_*z_1e^{t_*z_1} & e^{t_*z_2}+t_*z_2e^{t_*z_2} & \cdots & e^{t_*z_n}+t_*z_ne^{t_*z_n} \end{pmatrix}.
$
\end{center}
Its rank is strictly less than $2$ if and only if 
\begin{center}
$z_1=z_2=\cdots=z_n$.
\end{center}
Since $a_1\neq \cfrac{a_0-\log n}{t_*}$, this is excluded from points in $\Gamma_{\ref{eq35},\ref{eq36}}$. By the implicit function  theorem \cite{Lan}, if $\Gamma_{\ref{eq35},\ref{eq36}}^n$ is not empty, it is an $(n-2)$-dimension smooth manifold locally. The gradient of the function $Q_0(z_1,\cdots,z_n)-e^{a_0}$ is 
\begin{center}
$\bigtriangledown(Q_0(z_1,\cdots,z_n)-e^{a_0})=(t_*e^{t_*z_1}, t_*e^{t_*z_2}, \cdots, t_*e^{t_*z_n})$,
\end{center}
whose individual components will always be strictly positive. The gradient of the function $Q_1(z_1,\cdots,z_n)-a_1e^{a_0}$ is 
\begin{center}
$\bigtriangledown(Q_1(z_1,\cdots,z_n)-a_1e^{a_0})=(e^{t_*z_1}+t_*z_1e^{t_*z_1}, e^{t_*z_2}+t_*z_2e^{t_*z_2}, \cdots, e^{t_*z_n}+t_*z_ne^{t_*z_n})$,
\end{center}
with the $i$-th individual component vanishes if and only if $z_i=-\cfrac{1}{t_*}$ for $1\leq i\leq n$. So $\Gamma_{\ref{eq35}}^n$ and $\Gamma_{\ref{eq36}}^n$ cannot be tangent to each other. Moreover, note that 
\begin{center}
$e^{t_*z_i}+t_*z_ie^{t_*z_i}>0$
\end{center}
if $z_i>-\cfrac{1}{t_*}$ while
\begin{center}
$e^{t_*z_i}+t_*z_ie^{t_*z_i}<0$
\end{center}
if $z_i<-\cfrac{1}{t_*}$ for any $1\leq i\leq n$. These force the intersection of zeros of the two functions $Q_0(z_1,\cdots,z_n)-e^{a_0}$ and $Q_1(z_1,\cdots,z_n)-a_1e^{a_0}$ to be connected, if the intersection is not empty. This implies $\Gamma_{\ref{eq35},\ref{eq36}}$ is a manifold globally in case of being nonempty. $\Gamma_{\ref{eq35},\ref{eq36}}^n$ is compact since it is a bounded set.
\end{proof}
Let 
\begin{center}
$\Gamma_{\ref{eq35},1,2,1}^3=\{(z_1, z_2, z_3)\in\mathbb{R}^3: z_1, z_2, z_3 \mbox{ satisfy } e^{z_1}+e^{z_2}+e^{z_3}=e^2\}$
\end{center}
and
\begin{center}
$\Gamma_{\ref{eq36},1,2,1}^3=\{(z_1, z_2, z_3)\in\mathbb{R}^3: z_1, z_2, z_3 \mbox{ satisfy } z_1e^{z_1}+z_2e^{z_2}+z_3e^{z_3}=e^2\}$
\end{center}
be the corresponding surfaces with $t_*=1, a_0=2, a_1=1$. Figure \ref{fig2} depicts parts of the two $2$-dimension surfaces, whose intersection will be a $1$-dimension smooth curve.

Equipped with all the above results, now we are ready to prove Theorem \ref{thm3}.

\begin{proof}[Proof of Theorem \ref{thm3}]
First, for the given $t_*>0$ and $(a_0, a_1)\in\mathbb{R}^2$ satisfying (\ref{eq32}), if $n$ is large enough, $\Gamma_{\ref{eq35},\ref{eq36}}^n$ is not empty according to Corollary \ref{cor4}. So $\Gamma_{\ref{eq35},\ref{eq36}}^n$ is a compact $(n-2)$-dimension smooth manifold for $n$ large enough. In the following we always assume $n$ is large enough.  Now let
\begin{center}
$m_{t_*,a_0,a_1,n}=\min\Big\{\cfrac{R_2(t_*,z_1,z_2,\cdots,z_n)}{e^{a_0}}-a_1^2:(z_1, z_2, \cdots, z_n)\in\Gamma_{\ref{eq35},\ref{eq36}}^n\Big\}$
\end{center}
while
\begin{equation}\label{eq48}
M_{t_*,a_0,a_1,n}=\max\Big\{\cfrac{R_2(t_*,z_1,z_2,\cdots,z_n)}{e^{a_0}}-a_1^2:(z_1, z_2, \cdots, z_n)\in\Gamma_{\ref{eq35},\ref{eq36}}^n\Big\}.
\end{equation}
For any $m_{t_*,a_0,a_1,n}\leq a_2\leq M_{t_*,a_0,a_1,n}$, since $\Gamma_{\ref{eq35},\ref{eq36}}^n$ is a smooth manifold, there exist $\{c_{i,n}\}_{i=1}^n\subset\mathbb{R}$, such that $(c_{1,n},c_{2,n},\cdots,c_{n,n})$ satisfies $(\ref{eq35}),(\ref{eq36})$ and  
\begin{equation}\label{eq40}
a_2=\cfrac{Q_2(t_*,c_{1,n},\cdots,c_{n,n})}{Q_0^2(t_*,c_{1,n},\cdots,c_{n,n})}=\cfrac{R_2(t_*,c_{1,n},\cdots,c_{n,n})}{e^{a_0}}-a_1^2
\end{equation}
simultaneously. Now let
\begin{center}
$\phi(x)=c_{x_0,n}$
\end{center}  
for $x=\cdots x_{-1}x_{0}x_{1}\cdots\in[x_0]$ on the full shift space $\Lambda_n^{\mathbb{Z}}$. It is a locally constant potential. According to (\ref{eq39}) and (\ref{eq40}), we have
\begin{equation}\label{eq42}
P''(t_*\phi)=\cfrac{Q_2(t_*,c_{1,n},\cdots,c_{n,n})}{Q_0^2(t_*,c_{1,n},\cdots,c_{n,n})}=a_2.
\end{equation}
Since $(c_{1,n},c_{2,n},\cdots,c_{n,n})$ satisfies $(\ref{eq35})$ and $(\ref{eq36})$, we have
\begin{equation}\label{eq43}
P(t_*\phi)=\cfrac{Q_2(t_*,c_{1,n},\cdots,c_{n,n})}{Q_0^2(t_*,c_{1,n},\cdots,c_{n,n})}=a_0
\end{equation}
while
\begin{equation}\label{eq44}
P'(t_*\phi)=\cfrac{Q_2(t_*,c_{1,n},\cdots,c_{n,n})}{Q_0^2(t_*,c_{1,n},\cdots,c_{n,n})}=a_1.
\end{equation}
Note that $P(t\phi)$ is analytic with respect to $t$ on $(\alpha,\infty)$ for any $\alpha>0$, so there exists some $\delta_n>0$, such that (\ref{eq41}) holds on $[t_*-\delta_n,t_*+\delta_n]$, considering (\ref{eq42}), (\ref{eq43}) and (\ref{eq44}).
\end{proof}

In the following we illustrate some dependent relationship between
\begin{center}
$\{m_{t_*,a_0,a_1,n}, M_{t_*,a_0,a_1,n}\}_{n\in\mathbb{N}}$ 
\end{center}
and some particular $t_*,a_0,a_1,n$ satisfying (\ref{eq32}). There should be some universal relationship between
them, while we hope the following observations will provide some hints. The first one is that it is possible for $m_{t_*,a_0,a_1,n}=0$ for some  $t_*,a_0,a_1,n$.

\begin{pro}
Let $t_*>0$ and $(a_0,a_1)\in\mathbb{R}^2$ satisfy (\ref{eq32}). Then $m_{t_*,a_0,a_1,n}=0$ for $n\geq 2$ if and only if 
\begin{equation}\label{eq45}
a_1=\cfrac{a_0-\log n}{t_*}.
\end{equation}
\end{pro}
\begin{proof}
Note that $m_{t_*,a_0,a_1,n}=0$ is equivalent to say that there exists some locally constant potential $\phi$ on $\Lambda_n^{\mathbb{Z}}$ such that $P''(t_*\phi)=0$ according to Theorem \ref{thm3}. By \cite[Proposition 4.12]{PP}, this happens if and only if $\phi$ is a constant potential on $\Lambda_n^{\mathbb{Z}}$. In this case we have 
\begin{center}
$\phi(x)=\cfrac{a_0-\log n}{t_*}$
\end{center}
for any $x\in\Lambda_n^{\mathbb{Z}}$, which implies (\ref{eq45}).
\end{proof}

This result does not tell things about the sequence 
\begin{center}
$\{m_{t_*,a_0,a_1,n}\}_{\ n\in\mathbb{N} \mbox{ large enough }}$ 
\end{center}
for given  $t_*,a_0,a_1$, since (\ref{eq45}) will never be true for any $n$ large enough for fixed $t_*,a_0,a_1$. The following result describes some limit behaviour of the sequence 
\begin{center}
$\{M_{t_*,a_0,a_1,n}\}_{\ n\in\mathbb{N} \mbox{ large enough }}$ 
\end{center}
for $t_*=1,a_0=2,a_1=1$.

\begin{pro}\label{pro2}
Let $t_*=1,a_0=2,a_1=1$, in symbols of Theorem \ref{thm3}, we have
\begin{equation}
\lim_{n\rightarrow\infty} M_{1,2,1,n}=\infty.
\end{equation}
\end{pro}

To justify Proposition \ref{pro2},  we first illustrate some basic properties about the function $ze^{t_*z}$ for $t_*>0$.
\begin{lemma}
For $t_*>0$, $ze^{t_*z}$ is strictly decreasing on $(-\infty,-\cfrac{1}{t_*})$, strictly increasing on $(-\cfrac{1}{t_*},\infty)$, while it attains its minimum $-\cfrac{1}{t_*}e^{-1}$ at $z=-\cfrac{1}{t_*}$. It admits one and only one inflection in $(-\infty,-\cfrac{1}{t_*})$.
\end{lemma}
\begin{proof}
One can check these conclusions by some direct computations on the first and second derivatives of the function $ze^{t_*z}$.
\end{proof}

\begin{figure}[h]
\centering
\includegraphics[scale=0.84]{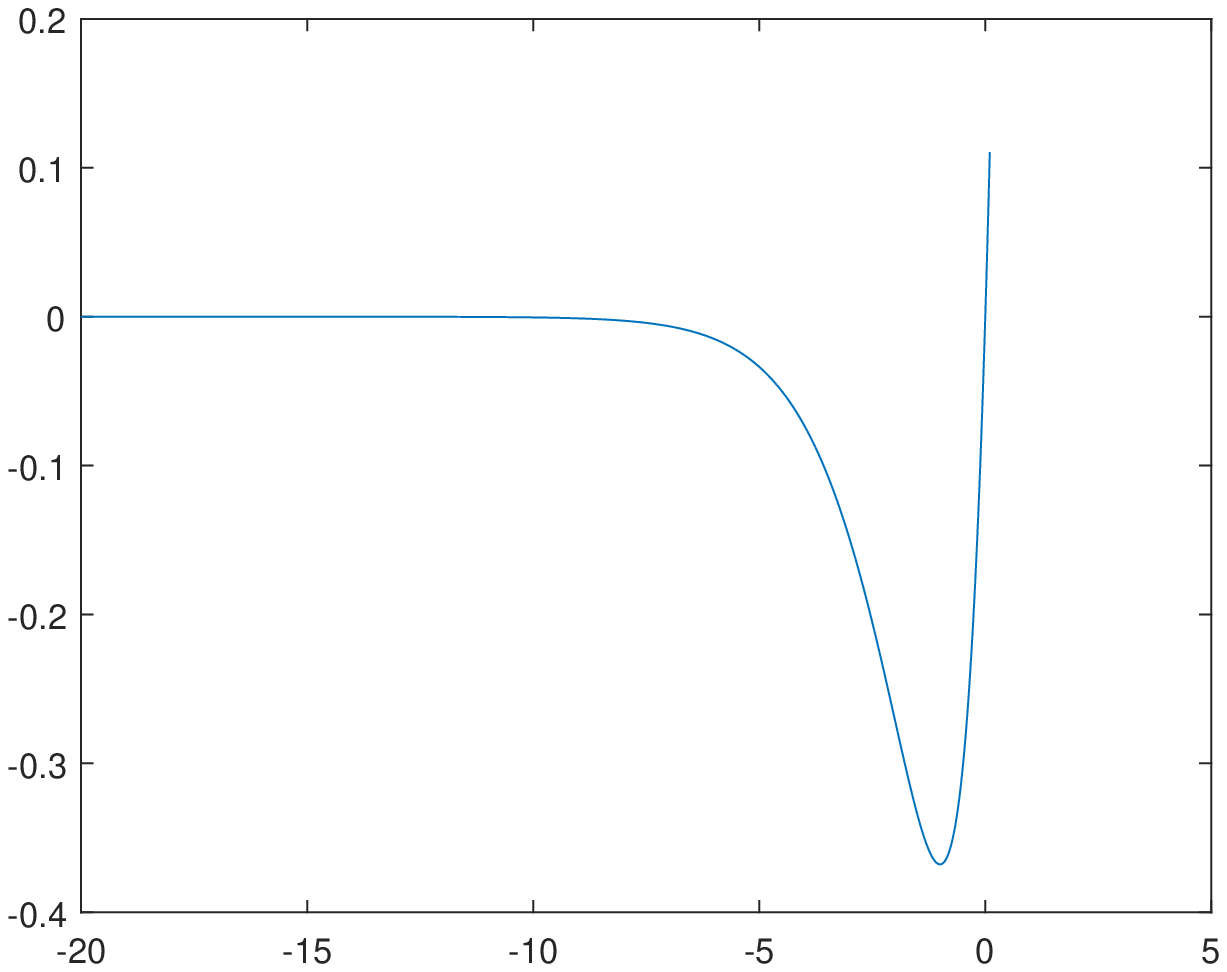}
\caption{Graph of $\varsigma(z)=ze^{z}$}
\label{fig3}
\end{figure}

In Figure \ref{fig3} we depict the graph of $\varsigma(z)=ze^{z}$.

\begin{figure}[h]
\centering
\includegraphics[scale=0.85]{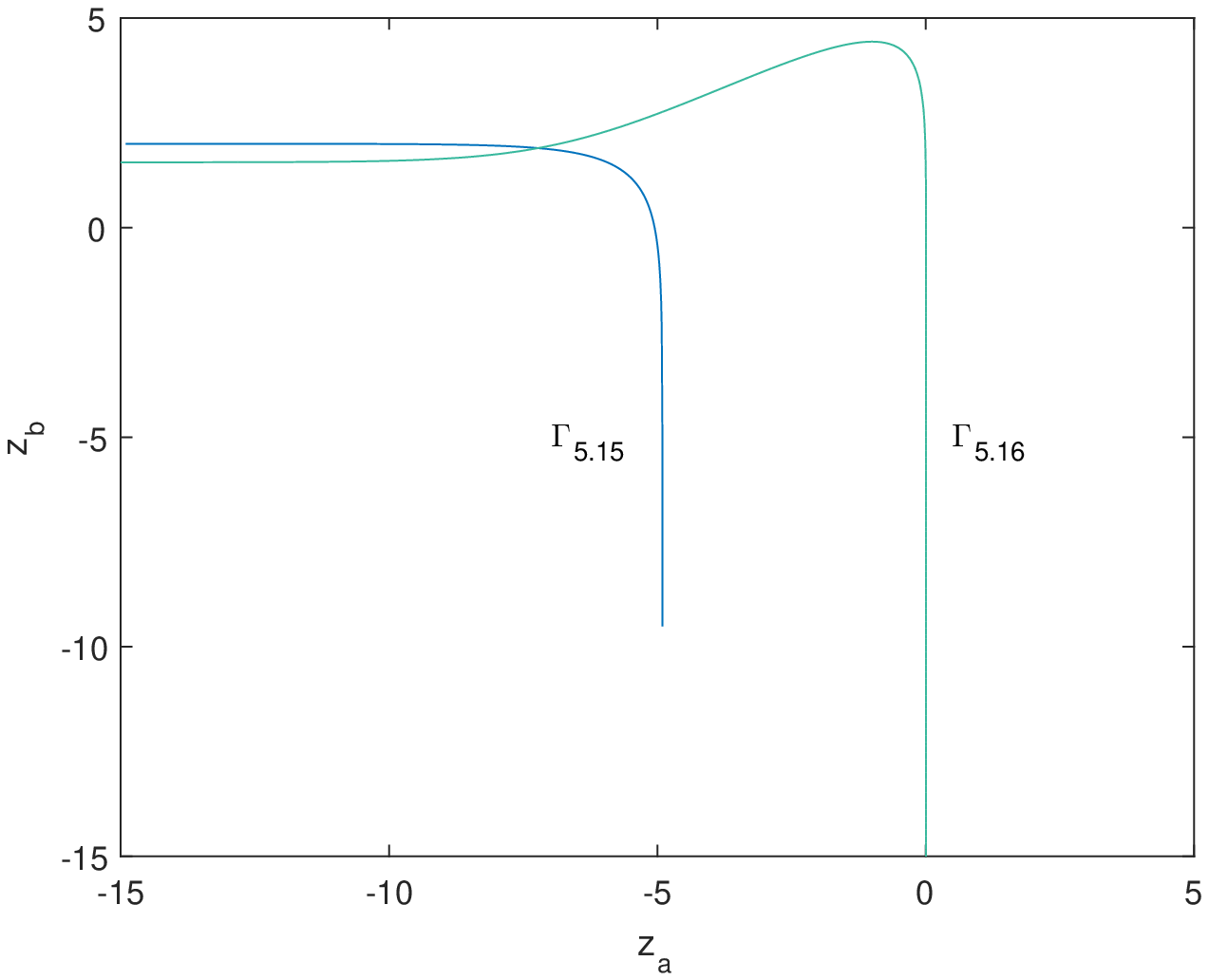}
\caption{$\Gamma_{\ref{eq46}}$ and $\Gamma_{\ref{eq47}}$}
\label{fig4}
\end{figure}

\begin{proof}[Proof of Proposition \ref{pro2}]
Since we are considering the limit behaviour of $M_{1,2,1,n}$, we always assume $n$ is large enough throughout the proof. Now consider the following two equations 
\begin{equation}\label{eq46}
(n-1)e^{z_a}+e^{z_b}=e^2
\end{equation} 
and
\begin{equation}\label{eq47}
(n-1)z_ae^{z_a}+z_be^{z_b}=e^2
\end{equation}
with unknowns $z_a, z_b$. Let 
\begin{center}
$\Gamma_{\ref{eq46}}=\{(z_a, z_b)\in\mathbb{R}^2: z_a, z_b \mbox{ satisfy } (\ref{eq46})\}$
\end{center}
and
\begin{center}
$\Gamma_{\ref{eq47}}=\{(z_a, z_b)\in\mathbb{R}^2: z_a, z_b \mbox{ satisfy } (\ref{eq47})\}$.
\end{center}
We describe the graph of $\Gamma_{\ref{eq46}}$ and $\Gamma_{\ref{eq47}}$ separately in the following. $\Gamma_{\ref{eq46}}$ is a $1$ dimensional smooth curve with two asymptotes $z_a=2-\log(n-1)$ and $z_b=2$. It is strictly decreasing  when we consider the curve as the graph of the function
\begin{center}
$z_b=\log\big(e^2-(n-1)e^{z_a}\big)$
\end{center}
for $z_a\in(-\infty, 2-\log(n-1))$. $\Gamma_{\ref{eq47}}$ is also a $1$ dimensional smooth curve with two asymptotes $z_a=\varsigma^{-1}(\cfrac{e^2}{n-1})$ and $z_b=\varsigma^{-1}(e^2)$. When we consider the $\Gamma_{\ref{eq47}}$ as the graph of the function
\begin{center}
$z_b=\eta(z_a)$
\end{center}
as the implicit function induced by $(\ref{eq47})$, it is strictly increasing for $z_a\in(-\infty,-1)$, strictly decreasing for $z_a\in(-1,\varsigma^{-1}(\cfrac{e^2}{n-1}))$, with its maximum  $\varsigma^{-1}(e^2+(n-1)e^{-1})$ attained at $z_a=-1$. Let $\varsigma_l^{-1}(-\cfrac{e^2}{n-1})$ be the smaller one of the two intersections of $z_b=2$ and $\Gamma_{\ref{eq47}}$, then $\Gamma_{\ref{eq46}}$ and $\Gamma_{\ref{eq47}}$ must intersection at some unique point $c_{a,n}\in(-\infty, \varsigma_l^{-1}(-\cfrac{e^2}{n-1}))$. Obviously 
\begin{center}
$\lim_{n\rightarrow\infty} c_{a,n}=-\infty$
\end{center} 
since $\lim_{n\rightarrow\infty} \varsigma_l^{-1}(-\cfrac{e^2}{n-1})=-\infty$. Now we analyse the order of $c_{a,n}$ with respect to $n$ as $n\rightarrow\infty$. Let 
\begin{center}
$z_{a,n}=-\log n-\log\log n+\log 1-1$.
\end{center}
One can check that
\begin{center}
$\lim_{n\rightarrow\infty, z_b\rightarrow 2} ((n-1)e^{z_{a,n}}+e^{z_b})=e^2$
\end{center}
while
\begin{center}
$\lim_{n\rightarrow\infty, z_b\rightarrow 2} ((n-1)z_{a,n}e^{z_{a,n}}+z_be^{z_b})=e^2$.
\end{center}
These imply that 
\begin{center}
$c_{a,n}=-\log n-\log\log n+o(\log\log n)$.
\end{center}
Note that  $(c_{a,n},c_{a,n},\cdots, c_{a,n}, \eta(c_{a,n}))\in \Gamma_{\ref{eq35},\ref{eq36}}^n$ for $t_*=1, a_0=2, a_1=1$.  Now 
\begin{center}
$
\begin{array}{ll}
& R_2(c_{a,n},c_{a,n},\cdots, c_{a,n}, \eta(c_{a,n}))\vspace{3mm}\\
= & (n-1)c_{a,n}^2e^{c_{a,n}}+(\eta(c_{a,n}))^2e^{\eta(c_{a,n})}\vspace{3mm}\\
= & (n-1)(-\log n-\log\log n+o(\log\log n))^2e^{-\log n-\log\log n+o(\log\log n)}+4e^2+o(1)\vspace{3mm}\\
= & \log n+o(\log n),
\end{array}
$
\end{center}
from which it is easy to see that 
\begin{center}
$\lim_{n\rightarrow\infty} R_2(c_{a,n},c_{a,n},\cdots, c_{a,n}, \eta(c_{a,n}))=\infty$.
\end{center}
This forces
\begin{center}
$\lim_{n\rightarrow\infty} M_{1,2,1,n}=\infty$,
\end{center}
considering (\ref{eq48}).

\end{proof}

We provide the readers with the curves $\Gamma_{\ref{eq46}}$ and $\Gamma_{\ref{eq47}}$ in Figure \ref{fig4}. Obviously some more general conclusions are available if one considers variations of the parameters  $t_*,a_0,a_1$ in Proposition \ref{pro2}. At last we provide the readers with some solutions $\{c_{a,n}\}_{n\in\mathbb{N}}$ and $\{\eta(c_{a,n})\}_{n\in\mathbb{N}}$ in Table \ref{tab3}, from which one can see the order of decay and increase of the sequences with respect to $n$ clearly.

\begin{table}[htbp]
\centering
\footnotesize
\caption{$\{c_{a,n}\}_{n\in\mathbb{N}}$ and $\{\eta(c_{a,n})\}_{n\in\mathbb{N}}$}
\label{tab3}
\renewcommand{\arraystretch}{1}
\begin{tabular}{|c|c|c|}
\hline
$n$ &  $c_{a,n}$ & $\eta(c_{a,n})$ \\
\hline
$10$     &   -1.8599539391797653780996686364493    &    1.7634042477581860636342812520981\\
\hline
$10^2$   &   -4.6278529940301947157458180305676    &    1.8580906928560505140960875180438\\
\hline
$10^3$   &   -7.2278923365046354303919671475052    &    1.8965708210067454817129699066334\\
\hline
$10^4$   &   -9.7529279223041958189401940128674    &    1.9180710389285259082138396366755\\
\hline 
$10^5$   &   -12.23426184122178540565187685582     &     1.9319494203818796717151866525306\\
\hline  
$10^6$   &   -14.686689485112383196253350885528     &     1.941701042038176132682488585943\\
\hline  
$10^7$   &   -17.118475509130338419321449219176     &     1.9489507180131363431129601417792\\
\hline 
$10^8$   &   -19.534737736752111249670741176574     &     1.9545628133690736391913141129777\\
\hline 
$10^9$   &   -21.938877884281897893422087428599     &     1.9590417833080193886068703580662\\
\hline  
$10^{10}$   &   -24.333277592346602338263750350022     &     1.9627027620469153955488959845337\\
\hline  
$10^{11}$   &   -26.719672172461371813735932628894     &     1.9657531814729595378854181456218\\
\hline   
$10^{12}$   &   -29.099366670257435261982274861811     &     1.9683353707111573738492465130807\\
\hline  
$10^{13}$   &   -31.473368167571030624456199153849     &     1.970550350496947761285545176838\\
\hline  
$10^{14}$   &   -33.842470627269595326611535858951     &     1.9724718685216929582206115029034\\
\hline  
$10^{15}$   &   -36.20731141238751139407393422892      &     1.9741550583546827046855344007126\\
\hline  
$10^{16}$   &   -38.568410155198951836337896822881      &    1.97564198636943790477268372057\\
\hline   
$10^{17}$   &   -40.926196222869058989174011616314      &    1.9769653208730088904749619599928\\
\hline  
$10^{18}$   &   -43.28102858421294787781225809291      &     1.9781508271703613365389080750692\\
\hline   
$10^{19}$   &   -45.633210475623427729647938869856      &     1.9792191056459012534062976747755\\
\hline   
$10^{20}$   &   -47.983000423353389741328990557576      &     1.9801868284846851379610473178804\\
\hline 
$10^{21}$   &   -50.330620660008332271820694306839      &     1.9810676363715292020369862557429\\
\hline  
$10^{22}$   &   -52.676263643082855194671803053742      &     1.9818727996772032079642260800619\\
\hline 
$10^{23}$   &   -55.020097168291592849066888176454      &     1.9826117134133018944596936081392\\
\hline   
$10^{24}$   &   -57.362268427077060922578379063246      &     1.9832922728467949817312209115653\\
\hline 
$10^{25}$   &   -59.702907260160132201351723856461      &     1.9839211621102961084222105523408\\
\hline  
$10^{26}$   &   -62.042128791447074538616865826092      &     1.9845040784885601043186175801529\\
\hline 
$10^{27}$   &   -64.380035579030470553978577616248      &     1.9850459085371711281404342988732\\
\hline  
$10^{28}$   &   -66.716719386002755126963619613768      &     1.9855508677057357884921072471682\\
\hline  
$10^{29}$   &   -69.052262649137714881922574449762      &     1.9860226120088820356292880321385\\
\hline 
$10^{30}$   &   -71.386739705385277326962820249044      &     1.9864643280735340181774784668139\\
\hline 
$10^{31}$   &   -73.7202178226698188293210417966        &     1.9868788063025456510398996161088\\
\hline 
$10^{32}$   &   -76.052758071376257724956806229201      &     1.9872685007417625212636588061233\\
\hline 
$10^{33}$   &   -78.384416065240707497606345034329      &     1.9876355783911370649894789906831\\
\hline 
$10^{34}$   &   -80.715242594490126828808238297291      &     1.9879819600725889180558042388717\\
\hline 
$10^{35}$   &   -83.045284169538297201269228695051      &     1.9883093544968926933418986392956\\
\hline 
$10^{36}$   &   -85.374583490011093204926910773042      &     1.9886192868162322855194994642595\\
\hline 
$10^{37}$   &   -87.703179851099500408441885242821      &     1.9889131226778662869710690898493\\
\hline 
$10^{38}$   &   -90.031109497045012553979249690171      &     1.9891920885858755212588911444123\\
\hline 
$10^{39}$   &   -92.358405929815521230622254914238      &     1.9894572892164831961499157326311\\
\hline 
$10^{40}$   &   -94.685100179630439886817678169988      &     1.9897097222064583485549741614556\\
\hline 
\end{tabular}
\end{table}

\end{document}